\documentclass[11pt,a4paper]{article}

\usepackage{epsf,epsfig,amsfonts,amsgen,amsmath,amstext,amsbsy,amsopn,amsthm
}
\usepackage{amsmath,times,mathptmx}
\usepackage{amsfonts,amsthm,amssymb}
\usepackage{amsfonts}
\usepackage{graphics}
\usepackage{latexsym,bm}
\usepackage{amsfonts,amsthm,amssymb,bbding}
\usepackage{indentfirst}
\usepackage{graphicx}
\usepackage{color}
\usepackage[colorlinks=true,anchorcolor=blue,filecolor=blue,linkcolor=blue,urlcolor=blue,citecolor=blue]{hyperref}
\usepackage{float}
\usepackage{tikz}
\newtheorem{thm}{Theorem}

\newtheorem{lemma}{Lemma}
\newtheorem{false statement}{False statement}

\theoremstyle{definition}

\newtheorem{claim}{Claim}

\newtheorem{case}{Case}

\begin{document}

\title{Extremal distance spectral radius of graphs with $h$-extra $r$-component connectivity
\footnote{Supported by Natural Science Foundation of Xinjiang Uygur Autonomous Region (No. 2024D01C41), NSFC (No. 12361071), the Basic Scientific Research in Universities of Xinjiang Uygur Autonomous Region (No. XJEDU2025P001), Tianshan Talent Training Program (No. 2024TSYCCX0013).}}
\author{{Daoxia Zhang, Dan Li$^a$}\thanks{Corresponding author. E-mail: ldxjedu@163.com.}, {Wenxiu Ding}\\
{\footnotesize $^a$College of Mathematics and System Science, Xinjiang University, Urumqi 830046, China}}
\date{}

\maketitle {\flushleft\large\bf Abstract:}
For two integers $r\geq 2$ and $h\geq 0$, the $h$-extra $r$-component connectivity of a graph $G$, denoted by $c\kappa_{r}^{h}$, is defined as the minimum number of vertices whose removal produces a disconnected graph with at least $r$ components, where each component contains at least $h+1$ vertices. Let $\mathcal{G}_{n,\delta}^{c\kappa_{r}^{h}}$ represent the set of graphs of order $n$ with minimum degree $\delta$ and $h$-extra $r$-component connectivity $c\kappa_{r}^{h}$. Hu, Lin, and Zhang [\textit{Discrete Math.} \textbf{345} (2025) 114621] investigated the case when $h=0$ within $\mathcal{G}_{n,\delta}^{c\kappa_{r}^{h}}$, and characterized the corresponding extremal graphs that minimize the distance spectral radius. In this paper, we further explore the relevant extremal graphs in $\mathcal{G}_{n,\delta}^{c\kappa_{r}^{h}}$ for $h\geq 1$. 

\vspace{0.1cm}
\begin{flushleft}
\end{flushleft}

\section{Introduction}
Throughout this paper, we only consider simple and undirected graphs. Let $G$ be a graph  with vertex set $V(G)=\{v_{1}, v_{2},\ldots,v_{n}\}$ and edge set $E(G)$. 
The distance between $v_{i}$ and $v_{j}$, denoted by $d_{G}(v_{i},v_{j})$ (or $d_{ij}$), is the length of a shortest path from $v_{i}$ to $v_{j}$. The distance matrix of $G$, denoted by $D(G)$, is the $n\times n$ real symmetric matrix whose $(i,j)$-entry is $d_{G}(v_{i},v_{j})$ (or $d_{ij}$), then we can order the eigenvalues of $D(G)$ as $\lambda_{1}(G)\geq \lambda_{2}(G)\geq \cdots \geq \lambda_{n}(G)$. By Perron-Frobenius theorem, $\lambda_{1}(G)$ is always positive (unless $G$ is trivial) and $\lambda_{1}(G)\geq |\lambda_{i}(G)|$ for $i=2, 3, \ldots, n$. We call $\lambda_{1}(G)$ the distance spectral radius. Denote by $\delta(G)$ the minimum degree of $G$.


Over the past decades, the distance spectral radius of graphs with given parameters inspired much interest and attracted the attention of researchers. The connectivity $\kappa(G)$ of a graph $G$ is the minimum number of vertices whose removal results in a disconnected graph or in the trivial graph. Connectivity and its generalizations have been studied due to their impact on the fault tolerance and diagnosability of interconnection networks {\cite{M.J. Ma, L.N. Zhao}}. In 2010, Liu \cite{Z.Z. Liu} characterized graphs with minimum distance spectral radius among $n$-vertex graphs with fixed vertex connectivity (resp. matching number, chromatic number). Zhang and Godsil \cite{X.L. Zhang} determined the graph with $k$ cut vertices with minimal distance spectral radius. In 2012, Nath and Paul \cite{M. Nath1} characterized the graphs with minimal distance spectral radius in the class of all connected bipartite graphs with a given vertex connectivity. Subsequently, they \cite{M. Nath} determined the graphs with minimum distance spectral radius among $n$-vertex graphs with given connectivity and minimum degree. Zhang \cite{Zhang} determined the $n$-vertex graphs of given diameter with the minimum distance spectral radius. Zhang, Li and Gutman \cite{M.J. Zhang} generalized this result by characterizing the graphs of order $n$ with given connectivity and diameter having minimum distance spectral radius. In addition, they determined the minimum distance spectral radius of graphs among the $n$-vertex graphs with given connectivity and independence number, and characterized the corresponding extremal graph. Recently, Hu and Zhang \cite{Hu} established a distance spectral condition for a graph to be $1$-binding and characterized the corresponding extremal graphs. For more results, one may refer to {\cite{D.D. Fan, L. Zhang, Y.K. Zhang}}.


In 1983, Harary \cite{F. Harary} introduced the concept of conditional connectivity by imposing some conditions on the components of $G-S$, where $S$ is a subset of edges or vertices. Recently, The extremal problem regarding the conditional connectivity of graphs has been well investigated by many researchers. In 2021, Li, Lan, Ning ,Tian, Zhang and Zhu \cite{B. Li} introduced the concept of $h$-extra $r$-component connectivity of a graph $G$, as an extension of the classic connectivity. For two integers $r\geq 2$ and $h\geq 0$, a subset $S$ is called an $h$-extra $r$-component cut of $G$ if there are at least $r$ connected components in $G-S$ and each component has at least $h+1$ vertices. The minimum size of any $h$-extra $r$-component cut of $G$, if exist, is the $h$-extra $r$-component connectivity of $G$, denoted as $c\kappa_{r}^{h}(G)$. Notice that $c\kappa_{2}^{0}(G)=\kappa(G)$. In 2024, Fan, Gu and Lin \cite{D.D. Fan(1)} took the lead in studying the combination of $r$-component connectivity (i.e., $h=0$) and spectral radius, and characterized the graphs that achieve the maximum spectral radius among all graphs of order $n$ with given minimum degree and $r$-component connectivity. For results combining conditional connectivity and spectral radius, we refer the readers to {\cite{W.X. Ding, Y. Wang(1), Y. Wang, Y. Wang(2)}}.


Inspired by the above conclusions, it is natural and interesting to further investigate the distance spectral condition corresponding to other conditional connectivity of graphs. Denote by $\mathcal{G}_{n,\delta}^{c\kappa_{r}^{h}}$ the set of graphs of order $n$ with minimum degree $\delta$ and $h$-extra $r$-component connectivity $c\kappa_{r}^{h}$. In 2025, Hu, Lin and Zhang \cite{Y.L. Hu} characterized the extremal graph that minimizes the distance spectral radius in $n$-vertex connected graphs with $h$-extra $r$-component connectivity. When $h=0$, they generalized the results of Fan, Gu and Lin \cite{D.D. Fan(1)} to the distance  spectrum. In this paper, we further investigate the corresponding extremal graphs in $\mathcal{G}_{n,\delta}^{c\kappa_{r}^{h}}$ for $h\geq 1$. As usual, $K_{n}$ denotes the complete graph of order $n$. For two vertex-disjoint graphs $G_{1}$ and $G_{2}$, we denote by $G_{1}\cup G_{2}$ the disjoint union of $G_{1}$ and $G_{2}$. The join $G_{1}\vee G_{2}$ is the graph obtained from $G_{1}\cup G_{2}$ by adding all possible edges between $V(G_{1})$ and $V(G_{2})$. It is a well-known fact that $\kappa(G)\geq \delta(G)$. Nevertheless, there is no bound between $\delta$ and $\kappa_{r}^{h}$ for $r\geq 3$. Let $G_{n,(h+1)^{r-1}}^{0, \delta}$ be the graph obtained from $K_{1}\cup (K_{c\kappa_{r}^{h}}\vee(K_{h}\cup (r-2)K_{h+1}\cup K_{n-c\kappa_{r}^{h}-(r-1)(h+1)}))$ by adding $\delta$ edges between $K_{1}$ and $K_{h}$. Let $G_{n,(h+1)^{r-1}}^{\delta-h, h}$ be the graph obtained from $K_{1}\cup(K_{c\kappa_{r}^{h}}\vee (K_{h}\cup(r-2) K_{h+1}\cup K_{n-c\kappa_{r}^{h}-(r-1)(h+1)}))$ by adding $h$ edges between $K_{1}$ and $K_{h}$, and then by adding $\delta-h$ edges between $K_{1}$ and $K_{c\kappa_{r}^{h}}$. Then, our result is as follows.

\begin{thm}\label{thm1.1}
\vspace*{2mm}
Suppose that $h\geq 1$ is an integer. Let  $G\in \mathcal{G}_{n,\delta}^{c\kappa_{r}^{h}}$ with $n\geq c\kappa_{r}^{h}+r(h+1)$. Then the following statements hold.\\[1.5ex]
(i) If $\delta \leq h$, then $\lambda_{1}(G) \geq \lambda_{1}(G_{n,(h+1)^{r-1}}^{0, \delta})$ with equality if and only if $G\cong G_{n,(h+1)^{r-1}}^{0, \delta}$;\\[1.5ex]
(ii) If $h< \delta< c\kappa_{r}^{h}+h$, then $\lambda_{1}(G) \geq \lambda_{1}(G_{n,(h+1)^{r-1}}^{\delta-h, h})$ with equality if and only if $G\cong G_{n,(h+1)^{r-1}}^{\delta-h, h}$;\\[1.5ex]
(iii) If $\delta \geq c\kappa_{r}^{h}+h$, then $\lambda_{1}(G) \geq \lambda_{1}(K_{c\kappa_{r}^{h}}\vee (K_{n-c\kappa_{r}^{h}-(r-1)(\delta-c\kappa_{r}^{h}+1)}\cup (r-1)K_{\delta-c\kappa_{r}^{h}+1})$ with equality if and only if $G\cong K_{c\kappa_{r}^{h}}\vee (K_{n-c\kappa_{r}^{h}-(r-1)(\delta-c\kappa_{r}^{h}+1)}\cup (r-1)K_{\delta-c\kappa_{r}^{h}+1}$).
\end{thm}

\section{Proof of Theorem \ref{thm1.1}}
\begin{lemma}\cite{C.D. Godsil}\label{lem1.1}
Let $e$ be an edge of $G$ such that $G-e$ is connected. Then $\lambda_{1}(G) < \lambda_{1}(G-e)$.
\end{lemma}

\begin{lemma}\cite{Y.L. Hu}\label{lem1.2}
Let $n$, $c$, $s$, $p$ and $n_{i}$ ($1\leq i\leq c$) be positive integers with $n=s+\sum_{i=1}^{c}n_{i}, n_{1}\geq n_{2} \geq \cdots \geq n_{c}\geq p$ and $n_{1}< n-s-(c-1)p$. Then $$\lambda_{1}(K_{s}\vee(K_{n_{1}}\cup K_{n_{2}}\cup \cdots \cup K_{n_{c}})> \lambda_{1}(K_{s}\vee(K_{n-s-p(c-1)}\cup(c-1)K_{p})).$$
\end{lemma}

\noindent{\bf{Proof of Theorem \ref{thm1.1}}.} Suppose that $G \in \mathcal{G}_{n,\delta}^{c\kappa_{r}^{h}}$ is a graph that attains the minimum distance spectral radius, where $r \geq 2$ and $h \geq 1$. Since $n \geq c\kappa_{r}^{h}+r(h+1)$, by the definition of $h$-extra $r$-component connectivity, there exists some nonempty subset $S \subseteq V(G)$ with $|S|=c\kappa_{r}^{h}$ such that $G-S$ contains at least $r$ components and each component has at least $h+1$ vertices. Let $G-S=B_{1} \cup B_{2} \cup \cdots \cup B_{q}$ where $q \geq r$, and let $|V(B_{i})|=n_{i}\geq 2$ for $1 \leq i \leq q$. Then $q=r$ by Lemma \ref{lem1.1}. Choose a vertex $u \in V(G)$ such that $d_{G}(u)=\delta$. Suppose that $d_{S}(u)=t$, $P_{i}=N_{G}(u) \cap V(B_{i})$ and $|P_{i}|=p_{i}$ for $1 \leq i \leq r$. Thus, $\sum_{i=1}^{r}p_{i}=\delta-t$. Without loss of generality, suppose that $n_{1}\leq n_{2} \leq \cdots \leq n_{r}$. Now, we divide the proof into the following three cases.

\begin{case}\label{case1.1}   
$\delta \leq h$.
\end{case}

Due to the uncertainty of minimum degree vertex $u$, we first present the following Claim.

\begin{claim}\label{claim1.1}
$u\notin S$.     
\end{claim} 
Otherwise, $u \in S$. By the minimality of $\lambda_{1}(G)$, $G-u\cong K_{c\kappa_{r}^{h}-1}\vee(K_{n_{1}}\cup K_{n_{2}}\cup \cdots \cup K_{n_{r}})$.  We first assert that $d_{S}(u)=t\neq \delta$. Otherwise, let $N_{S}(u)=\{v_{1}, v_{2},\ldots, v_{\delta}\}$. Choose a vertex $w\in V(K_{n_{1}})$. Let $G'=G-uv_{\delta}+uw$. Obviously, $G'\in \mathcal{G}_{n,\delta}^{c\kappa_{r}^{h}}$. Let $X$ be the Perron vector of $D(G')$ and $x(v)$ denote the entry of $X$ corresponding to the vertex $v\in V(G')$. By symmetry, set $x(v)=x_{d}$ for any $v\in V(K_{\delta-1})$, $x(v)=x_{k}$ for any $v\in V(K_{c\kappa_{r}^{h}-\delta})$, $x(w)=x_{w}$ for the vertex $w$, $x(v)=x_{1}$ for any $v\in V(K_{n_{1}-1})$, $x(v)=x_{i}$ for any $v\in V(K_{n_{i}})$ where $2 \leq i \leq r$ and $x(u)=x_{u}$ for the vertex $u$. Then 
$$X=(\underbrace{x_{d},\dots,x_{d}}_{\delta -1},\underbrace{x_{k},\dots,x_{k}}_{c\kappa_{r}^{h}-\delta},x_{w},\underbrace{x_{1},\dots,x_{1}}_{n_{1}-1},\ldots,\underbrace{x_{i},\dots,x_{i}}_{n_{i}},\ldots,x_{u})^{T}.$$
By $D(G')X=\lambda_{1}(G')X$, we have
$$\lambda_{1}(G')x_{w}=x_{u}+(\delta-1)x_{d}+(c\kappa_{r}^{h}-\delta)x_{k}+(n_{1}-1)x_{1}+2\sum_{i=2}^{r}n_{i}x_{i},$$
$$\lambda_{1}(G')x_{k}=2x_{u}+(\delta-1)x_{d}+(c\kappa_{r}^{h}-\delta-1)x_{k}+x_{w}+(n_{1}-1)x_{1}+\sum_{i=2}^{r}n_{i}x_{i},$$
$$\lambda_{1}(G')x_{i}=2x_{u}+(\delta-1)x_{d}+(c\kappa_{r}^{h}-\delta)x_{k}+2x_{w}+2(n_{1}-1)x_{1}+2\sum_{j=2,j\neq i}^{r}n_{j}x_{j}+(n_{i}-1)x_{i},$$
$$\lambda_{1}(G')x_{u}=(\delta-1)x_{d}+2(c\kappa_{r}^{h}-\delta)x_{k}+x_{w}+2(n_{1}-1)x_{1}+2\sum_{i=2}^{r}n_{i}x_{i},$$
from which we get
\begin{equation*}
\begin{aligned}
&\quad (\lambda_{1}(G')+3)\Big(\sum_{i=2}^{r}n_{i}x_{i}-x_{u}\Big)\\
&=\Big(2\sum_{i=2}^{r}n_{i}-3\Big)x_{u}+\Big(\sum_{i=2}^{r}n_{i}-1\Big)(\delta-1)x_{d}+\Big(\sum_{i=2}^{r}n_{i}-2\Big)(c\kappa_{r}^{h}-\delta)x_{k}+\Big(2\sum_{i=2}^{r}n_{i}-1\Big)x_{w}\\
&\quad +2\Big(\sum_{i=2}^{r}n_{i}-1\Big)x_{1}+\sum_{i=2}^{r}\sum_{j=2,j\neq i}^{r}n_{i}n_{j}x_{j}+\sum_{i=2}^{r}n_{i}^{2}x_{i}\\
& >0
\end{aligned}
\end{equation*}
and
$$x_{w}-x_{k}=\frac{\sum_{i=2}^{r}n_{i}x_{i}-x_{u}}{\lambda_{1}(G')+1}.$$
Note that $D(G)-D(G')$ is\\
\begin{center}
$\begin{matrix}
\delta-1\\
1\\
c\kappa_{r}^{h}-\delta-1\\
1\\
n_{1}-1\\
n_{2}\\
\vdots\\
n_{r}\\
u
\end{matrix}
\begin{bmatrix}
  O&  O&  O&  O&  O&  O&  \cdots&  O&  O \\
  O&  0&  O&  0&  O&  O&  \cdots&  O&  -1 \\
  O&  O&  O&  O&  O&  O&  \cdots&  O&  O \\
  O&  0&  O&  0&  O&  O&  \cdots&  O&  1 \\
  O&  O&  O&  O&  O&  O&  \cdots&  O&  O \\
  O&  O&  O&  O&  O&  O&  \cdots&  O&  O \\
  \vdots&  \vdots&  \vdots&  \vdots&  \vdots&  \vdots&  \ddots&  \vdots&  \vdots \\
  O&  O&  O&  O&  O&  O&  \cdots&  O&  O \\
  O&  -1&  O&  1&  O&  O&  \cdots&  O&  0
\end{bmatrix}$. 
\end{center}
Then
\begin{equation*}
\begin{aligned}
\lambda_{1}(G)-\lambda_{1}(G')&\ge X^{T} (D(G)-D(G'))X\\
& =2x_{u}(x_{w}-x_{k})\\
& >0~~(\text{since}~x_{w}>x_{k}),
\end{aligned}
\end{equation*}
a contradiction. This demonstrates that $d_{S}(u)=t\neq \delta$. We next assert that $N_{G}(u)\nsubseteq V(S\cup B_{i})$ and $N_{G}(u)\nsubseteq V(B_{i})$ for $1 \leq i \leq r$. Otherwise, either $N_{G}(u)\subseteq V(S\cup B_{i})$ or $N_{G}(u)\subseteq V(B_{i})$ for $1 \leq i \leq r$. Then we can obtain a new set $S^{*}=S-\{u\}$ with $c\kappa_{r}^{h}=|S^{*}|=|S|-1$, which contradicts $c\kappa_{r}^{h}=|S|$. This implies that $d_{S}(u)=t\neq \delta-1$. If $d_{S}(u)=t\leq \delta-2$, without loss of generality, suppose that $p_{1}\geq 1$. (There must exist two integers $a,b \in [1,r]$ such that $N_{G}(u)\cap B_{a}\neq \emptyset$ and $N_{G}(u)\cap B_{b}\neq \emptyset$. But for the sake of conciseness in writing, suppose that $p_{1}\geq 1$.) By the minimality of $\lambda_{1}(G)$, we have $G-u\cong K_{c\kappa_{r}^{h}-1}\vee (K_{n_{1}}\cup K_{n_{2}}\cup \cdots \cup K_{n_{r}})$. Suppose that $P_{1}=\{w_{1},\ldots,w_{p_{1}}\}$ and $V(K_{n_{1}})\backslash P_{1}=\{w_{p_{1}+1},\ldots,w_{n_{1}}\}$. Let $G_{1}=G-\{uz \mid z\in P_{i}, 2\leq i\leq r\}+\{uw_{l} \mid p_{1}+1\leq l\leq \sum_{i=1}^{r}p_{i}\}+\{w_{n_{1}}w \mid w\in V(K_{n_{i}}), 2\leq i\leq r\}$. Obviously, $G_{1} \in \mathcal{G}_{n,\delta}^{c\kappa_{r}^{h}}$. Let $Y$ be the Perron vector of $D(G_{1})$ and $y(v)$ denote the entry of $Y$ corresponding to the vertex $v\in V(G_{1})$. By symmetry, set $y(v)=y_{d}$ for any $v\in V(K_{t})$, $y(v)=y_{k}$ for any $v\in V(K_{c\kappa_{r}^{h}-t})$, $y(v)=y_{1}'$ for any $v\in V(K_{\sum_{i=1}^{r}p_{i}})$, $y(v)=y_{1}$ for any $v\in V(K_{n_{1}-\sum_{i=1}^{r}p_{i}-1})$, $y(v)=y_{i}$ for any $v\in V(K_{n_{i}})$ where $2\leq i\leq r$ and $y(u)=y_{u}$ for the vertex $u$. Then
$$Y=(\underbrace{y_{d},\dots,y_{d}}_{t},\underbrace{y_{k},\dots,y_{k}}_{c\kappa_{r}^{h}-t},\underbrace{y_{1}',\dots,y_{1}'}_{\sum_{i=1}^{r}p_{i}},\underbrace{y_{1},\dots,y_{1}}_{n_{1}-\sum_{i=1}^{r}p_{i}-1},\ldots,\underbrace{y_{i},\dots,y_{i}}_{n_{i}},\ldots,y_{u})^{T}.$$
By $D(G_{1})Y=\lambda_{1}(G_{1})Y$, we have

$$\lambda_{1}(G_{1})y_{u}=ty_{d}+2(c\kappa_{r}^{h}-t)y_{k}+\sum_{i=1}^{r}p_{i}y_{1}'+2\Big(n_{1}-1-\sum_{i=1}^{r}p_{i}\Big)y_{1}+2\sum_{i=2}^{r}n_{i}y_{i},$$
$$\lambda_{1}(G_{1})y_{k}=ty_{d}+(c\kappa_{r}^{h}-t-1)y_{k}+2y_{u}+\sum_{i=1}^{r}p_{i}y_{1}'+\Big(n_{1}-1-\sum_{i=1}^{r}p_{i}\Big)y_{1}+\sum_{i=2}^{r}n_{i}y_{i},$$
$$\lambda_{1}(G_{1})y_{1}'=ty_{d}+(c\kappa_{r}^{h}-t)y_{k}+y_{u}+\Big(\sum_{i=1}^{r}p_{i}-1\Big)y_{1}'+\Big(n_{1}-1-\sum_{i=1}^{r}p_{i}\Big)y_{1}+2\sum_{i=2}^{r}n_{i}y_{i},$$
$$\lambda_{1}(G_{1})y_{i}=ty_{d}+(c\kappa_{r}^{h}-t)y_{k}+2y_{u}+2\sum_{i=1}^{r}p_{i}y_{1}'+2\Big(n_{1}-1-\sum_{i=1}^{r}p_{i}\Big)y_{1}+2\sum_{j=2,j\neq i}^{r}n_{j}y_{j}+(n_{i}-1)y_{i}.$$
Since $c\kappa_{r}^{h}> t$ and $n_{1}-1\geq h\geq \delta \geq \sum_{i=1}^{r}p_{i}$, combining the above equalities, we have
\begin{equation} \label{equ1}
\begin{aligned}
y_{u}-y_{k}=\frac{(c\kappa_{r}^{h}-t-1)y_{k}+(n_{1}-1-\sum_{i=1}^{r}p_{i})y_{1}+\sum_{i=2}^{r}n_{i}y_{i}}{\lambda_{1}(G_{1})+2}>0,
\end{aligned}
\end{equation}
\begin{equation} \label{equ2}
\begin{aligned}
y_{1}'-(y_{u}-y_{k})=\frac{ty_{d}+2y_{u}+\sum_{i=1}^{r}p_{i}y_{1}'+\sum_{i=2}^{r}n_{i}y_{i}}{\lambda_{1}(G_{1})+1}>0
\end{aligned}
\end{equation}
and
\begin{equation} \label{equ3}
\begin{aligned}
&\quad (\lambda_{1}(G_{1})+n_{i}+1)(y_{u}-y_{i})\\
& =(c\kappa_{r}^{h}-t)y_{k}+(n_{i}-1)y_{u}-\sum_{i=1}^{r}p_{i}y_{1}'\\
& =(c\kappa_{r}^{h}-t)y_{k}+(n_{i}-1)\left(y_{1}'+\frac{(c\kappa_{r}^{h}-t)y_{k}+(n_{1}-1-\sum_{i=1}^{r}p_{i})y_{1}}{\lambda_{1}(G_{1})+1}\right)-\sum_{i=1}^{r}p_{i}y_{1}'\\
& >(c\kappa_{r}^{h}-t)y_{k}+(n_{i}-1)y_{1}'-\sum_{i=1}^{r}p_{i}y_{1}'\\
& =(c\kappa_{r}^{h}-t)y_{k}+\Big(n_{i}-1-\sum_{i=1}^{r}p_{i}\Big)y_{1}'\\
& >0.
\end{aligned}
\end{equation}
Note that $D(G)-D(G_{1})$ is\\
\addtocounter{MaxMatrixCols}{10}
\begin{center}
$\begin{matrix}
t\\
c\kappa_{r}^{h}-t-1\\
1\\
p_{1}\\
\sum_{i=2}^{r}p_{i}\\
n_{1}-1-\sum_{i=1}^{r}p_{i}\\
p_{2}\\
n_{2}-p_{2}\\
\vdots\\
p_{r}\\
n_{r}-p_{r}\\
1
\end{matrix}
\begin{bmatrix}
  O&  O&  O&  O&  O&  O&  O&  O&  \cdots&  O&  O&  O \\
  O&  O&  O&  O&  O&  O&  O&  O&  \cdots&  O&  O&  O \\
  O&  O&  0&  O&  O&  O&  J&  J&  \cdots&  J&  J&  0 \\
  O&  O&  O&  O&  O&  O&  O&  O&  \cdots&  O&  O&  O \\
  O&  O&  O&  O&  O&  O&  O&  O&  \cdots&  O&  O&  J \\
  O&  O&  O&  O&  O&  O&  O&  O&  \cdots&  O&  O&  O \\
  O&  O&  J&  O&  O&  O&  O&  O&  \cdots&  O&  O& -J \\
  O&  O&  J&  O&  O&  O&  O&  O&  \cdots&  O&  O&  O \\
  \vdots&  \vdots&  \vdots&  \vdots&  \vdots&  \vdots&  \vdots&  \vdots&  \ddots&  \vdots&  \vdots&  \vdots \\
  O&  O&  J&  O&  O&  O&  O&  O&  \cdots&  O&  O&  -J \\
  O&  O&  J&  O&  O&  O&  O&  O&  \cdots&  O&  O&  O \\
  O&  O&  0&  O&  J&  O& -J&  O&  \cdots& -J&  O&  0
\end{bmatrix}$.
\end{center}
Combining this with (\ref{equ1})-(\ref{equ3}), we obtain
\begin{equation*}
\begin{aligned}
\lambda_{1}(G)-\lambda_{1}(G_{1})&\ge Y^{T} (D(G)-D(G_{1}))Y\\
& =2\Big(\sum_{i=2}^{r}n_{i}y_{i}y_{k}+\sum_{i=2}^{r}p_{i}y_{1}'y_{u}-\sum_{i=2}^{r}p_{i}y_{i}y_{u}\Big)\\
& >2\sum_{i=2}^{r}p_{i}\big(y_{1}'y_{u}-y_{i}(y_{u}-y_{k})\big)\\
& >2\sum_{i=2}^{r}p_{i}(y_{u}-y_{k})(y_{u}-y_{i})~~(\text{since}~y_{1}'>y_{u}-y_{k})\\
& >0~~(\text{since}~y_{u}>y_{k}~\text{and}~y_{u}>y_{i}).
\end{aligned}
\end{equation*}
Thus, $\lambda_{1}(G_{1})<\lambda_{1}(G)$, which contradicts the minimality of $\lambda_{1}(G)$. This implies that $u\notin S$.

\vspace{3mm}
By Claim \ref{claim1.1}, $u\notin S$, and so $u\in V(B_{i})$ for some $1\leq i\leq r$. Without loss of generality, we assume $u\in V(B_{1})$ for simplicity in subsequent notation an computation. Then $d_{B_{1}}(u)\geq 1$ because $n_{1}\geq h+1$ and $B_{1}$ is connected. 

\begin{claim}
$t=0$.     
\end{claim}

If not, assume that $d_{S}(u)=t\geq 1$. By the minimality of $\lambda_{1}(G)$, it follows that $G-u\cong K_{c\kappa_{r}^{h}}\vee(K_{n_{1}-1}\cup K_{n_{2}}\cup \cdots \cup K_{n_{r}})$. Suppose that
$P_{1}=\{w_{1},\ldots, w_{\delta-t}\}$ and $V(K_{n_{1}-1})\backslash P_{1}=\{w_{\delta-t+1},\ldots , w_{n_{1}-1}\}$. Let $G_{2}=G-\{uv\mid v\in N_{S}(u)\}+\{uw_{i}\mid \delta-t+1\leq i\leq \delta\}$. Obviously, $G_{2}\in \mathcal{G}_{n,\delta}^{c\kappa_{r}^{h}}$. Let $X$ be the Perron vector of $D(G_{2})$ and $x(v)$ denote the entry of $X$ corresponding to the vertex $v\in V(G_{2})$. By symmetry, set $x(v)=x_{k}$ for any $v\in V(K_{c\kappa_{r}^{h}})$, $x(v)=x_{1}'$ for any $v\in V(K_{\delta})$, $x(v)=x_{1}$ for any $v\in V(K_{n_{1}-\delta-1})$, $x(v)=x_{i}$ for any $v\in V(K_{n_{i}})$ where $2 \leq i \leq r$ and $x(u)=x_{u}$ for the vertex $u$. Then 
$$X=(\underbrace{x_{k},\dots,x_{k}}_{c\kappa_{r}^{h}},\underbrace{x_{1}',\dots,x_{1}'}_{\delta},\underbrace{x_{1},\dots,x_{1}}_{n_{1}-\delta-1},\dots,\underbrace{x_{i},\dots,x_{i}}_{n_{i}},\dots,x_{u})^{T}.$$
By $D(G_{2})X=\lambda_{1}(G_{2})X$, we have
$$\lambda_{1}(G_{2})x_{1}'=c\kappa_{r}^{h}x_{k}+x_{u}+(\delta-1)x_{1}'+(n_{1}-\delta-1)x_{1}+2\sum_{i=2}^{r}n_{i}x_{i},$$
$$\lambda_{1}(G_{2})x_{k}=(c\kappa_{r}^{h}-1)x_{k}+2x_{u}+\delta x_{1}'+(n_{1}-\delta-1)x_{1}+\sum_{i=2}^{r}n_{i}x_{i},$$
$$\lambda_{1}(G_{2})x_{i}=c\kappa_{r}^{h}x_{k}+2x_{u}+2\delta x_{1}'+2(n_{1}-\delta-1)x_{1}+2\sum_{j=2,j\neq i}^{r}n_{j}x_{j}+(n_{i}-1)x_{i},$$
$$\lambda_{1}(G_{2})x_{u}=2c\kappa_{r}^{h}x_{k}+\delta x_{1}'+2(n_{1}-\delta-1)x_{1}+2\sum_{i=2}^{r}n_{i}x_{i}.$$
Then
\begin{equation*}
\begin{aligned}
& \quad (\lambda_{1}(G_{2})+3)\Big(\sum_{i=2}^{r}n_{i}x_{i}-x_{u}\Big)\\
& =\Big(\sum_{i=2}^{r}n_{i}-2\Big)c\kappa_{r}^{h}x_{k}+\Big(2\sum_{i=2}^{r}n_{i}-3\Big)x_{u}+\Big(2\sum_{i=2}^{r}n_{i}-1\Big)\delta x_{1}'+2\Big(\sum_{i=2}^{r}n_{i}-1\Big)(n_{1}-\delta-1)x_{1}+\sum_{i=2}^{r}n_{i}^{2}x_{i}\\
& \quad +2\sum_{i=2}^{r}\sum_{j=2,j\neq i}^{r}n_{i}n_{i}x_{j}\\
& >0
\end{aligned}
\end{equation*}
and
$$x_{1}'-x_{k}=\frac{\sum_{i=2}^{r}n_{i}x_{i}-x_{u}}{\lambda_{1}(G_{2})+1}>0.$$
Note that $D(G)-D(G_{2})$ is\\
\begin{center}
$\begin{matrix}
t\\
c\kappa_{r}^{h}-t\\
\delta-t\\
t\\
n_{1}-\delta-1\\
n_{2}\\
\vdots\\
n_{r}\\
1
\end{matrix}
\begin{bmatrix}
  O&  O&  O&  O&  O&  O&  \cdots&  O&  -J \\
  O&  O&  O&  O&  O&  O&  \cdots&  O&  O \\
  O&  O&  O&  O&  O&  O&  \cdots&  O&  O \\
  O&  O&  O&  O&  O&  O&  \cdots&  O&  J \\
  O&  O&  O&  O&  O&  O&  \cdots&  O&  O \\
  O&  O&  O&  O&  O&  O&  \cdots&  O&  O \\
  \vdots&  \vdots&  \vdots&  \vdots&  \vdots&  \vdots&  \ddots&  \vdots&  \vdots \\
  O&  O&  O&  O&  O&  O&  \cdots&  O&  O \\
  -J&  O&  O&  J&  O&  O&  \cdots&  O&  0
\end{bmatrix}$.
\end{center}
Thus,
\begin{equation*}
\begin{aligned}
\lambda_{1}(G)-\lambda_{1}(G_{2})&\ge X^{T} (D(G)-D(G_{2}))X\\
&=2(tx_{1}'x_{u}-tx_{k}x_{u})\\
&=2tx_{u}(x_{1}'-x_{k})\\
&>0~~(\text{since}~x_{1}'>x_{k}),
\end{aligned}
\end{equation*}
a contradiction. Then we conclude that $t=0$.

Recall that $|V(B_{i})|=n_{i}\geq 2$ for $1\leq i\leq r$ and $n_{1}\leq n_{2}\leq \cdots \leq n_{r}$.

\begin{claim} \label{claim1.3}
$u\in V(B_{i})$ for some $n_{i}=h+1$, where $1\leq i\leq r$.     
\end{claim}

If $n_{r}=h+1$, the result is trivial. Next, we consider the case of $n_{r}>h+1$. We use proof by contradiction, assume that there exists $B_{q}$ such that $u\in V(B_{q})$ for $n_{q}>h+1$, where $1\leq q\leq r$. Then $G-u\cong K_{c\kappa_{r}^{h}}\vee (K_{n_{1}}\cup \cdots \cup K_{n_{q}-1}\cup \cdots \cup K_{n_{r}})$. Suppose that $P_{q}=\{w_{1},\ldots, w_{\delta}\}$ and $V(K_{n_{q}-1})\backslash P_{q}=\{w_{\delta+1},\ldots, w_{n_{q}-1}\}$. Let $G_{3}=G-\{w_{a}w_{b} \mid 1\leq a\leq h, h+1\leq b\leq n_{q}-1\}+\{w_{b}z \mid h+1\leq b\leq n_{q}-1, z\in V(K_{n_{j}})\}$, where $q\neq j$. Obviously, $G_{3}\in \mathcal{G}_{n,\delta}^{c\kappa_{r}^{h}}$. Let $Y$ be the Perron vector of $D(G_{3})$ and $y(v)$ denote the entry of $Y$ corresponding to the vertex $v\in V(G_{3})$. By symmetry, set $y(v)=y_{k}$ for any $v\in V(K_{c\kappa_{r}^{h}})$, $y(v)=y_{q}'$ for any $v\in V(K_{\delta})$, $y(v)=y_{q}$ for any $v\in V(K_{h-\delta})$, $y(v)=y_{j}$ for any $v\in V(K_{n_{q}+n_{j}-h-1})$, $y(v)=y_{l}$ for any $v\in V(K_{n_{l}})$ where $1 \leq l \leq r$ and $l\neq q,j$ and $y(u)=y_{u}$ for the vertex $u$. Then
$$Y=(\underbrace{y_{k},\dots,y_{k}}_{c\kappa_{r}^{h}},\underbrace{y_{1},\dots,y_{1}}_{n_{1}},\ldots,\underbrace{y_{q}',\dots,y_{q}'}_{\delta},\underbrace{y_{q},\dots,y_{q}}_{h-\delta},\ldots,\underbrace{y_{j},\dots,y_{j}}_{n_{q}+n_{j}-h-1},\ldots,\underbrace{y_{r},\dots,y_{r}}_{n_{r}},y_{u})^{T}.$$
By $D(G_{3})Y=\lambda_{1}(G_{3})Y$, we have
$$\lambda_{1}(G_{3})y_{j}=c\kappa_{r}^{h}y_{k}+2\sum_{l=1,l\neq q,j}^{r}n_{l}y_{l}+2y_{u}+2\delta y_{q}'+2(h-\delta)y_{q}+(n_{i}+n_{j}-h-2)y_{j},$$
$$\lambda_{1}(G_{3})y_{q}'=c\kappa_{r}^{h}y_{k}+2\sum_{l=1,l\neq q,j}^{r}n_{l}y_{l}+y_{u}+(\delta-1)y_{q}'+(h-\delta)y_{q}+2(n_{i}+n_{j}-h-2)y_{j},$$
$$\lambda_{1}(G_{3})y_{q}=c\kappa_{r}^{h}y_{k}+2\sum_{l=1,l\neq q,j}^{r}n_{l}y_{l}+2y_{u}+\delta y_{q}'+(h-\delta-1)y_{q}+2(n_{i}+n_{j}-h-2)y_{j}.$$
Then
\begin{equation*}
\begin{aligned}
&\quad (\lambda_{1}(G_{3})+n_{j}+1)\big(n_{j}y_{j}-\delta y_{q}'-(h-\delta)y_{q}\big)\\
& =(n_{j}-h)c\kappa_{r}^{h}y_{k}+2(n_{j}-h)\sum_{l=1,l\neq q,j}^{r}n_{l}y_{l}+(2n_{j}+\delta-2h)y_{u}+(n_{j}-h)\delta y_{q}'+(n_{j}-h)(h-\delta)y_{q}\\
& \quad +\big((n_{j}-2h)(n_{q}+n_{j}-h-1)+n_{j}^{2}\big)y_{j}\\
& >(n_{j}-h)c\kappa_{r}^{h}y_{k}+2(n_{j}-h)\sum_{l=1,l\neq q,j}^{r}n_{l}y_{l}+(2n_{j}+\delta-2h)y_{u}+(n_{j}-h)\delta y_{q}'+(n_{j}-h)(h-\delta)y_{q}\\
& \quad +2n_{j}(n_{j}-h)y_{j}\\
& >0.
\end{aligned}
\end{equation*}
Note that $D(G)-D(G_{3})$ is\\
\addtocounter{MaxMatrixCols}{10}
\begin{center}
$\begin{matrix}
c\kappa_{r}^{h}\\
n_{1}\\
\vdots\\
\delta\\
h-\delta\\
\vdots\\
n_{q}-h-1\\
n_{j}\\
\vdots\\
n_{r}\\
1
\end{matrix}
\begin{bmatrix}
  O&  O&  \cdots&   O&  O&    \cdots&  O&  O&  \cdots&  O&  O \\
  O&  O&  \cdots&    O&  O&   \cdots&  O&  O&  \cdots&  O&  O  \\
  \vdots&  \vdots&  \ddots&    \vdots&    \vdots&  \ddots&  \vdots&  \vdots&  \ddots&  \vdots&  \vdots \\
  O&  O&  \cdots&    O&  O&    \cdots& -J&  O&  \cdots&  O&  O \\
  O&  O&  \cdots&    O&  O&   \cdots& -J&  O&  \cdots&  O&  O \\
  \vdots&  \vdots&  \ddots&    \vdots&    \vdots&  \ddots&  \vdots&  \vdots&  \ddots&  \vdots&  \vdots \\
  O&  O&  \cdots&   -J& -J&    \cdots&  O&  J&  \cdots&  O&  O \\
  O&  O&  \cdots&    O&  O&    \cdots&  J&  O&  \cdots&  O&  O \\
  \vdots&  \vdots&  \ddots&   \vdots&   \vdots&  \ddots&  \vdots&  \vdots&  \ddots&  \vdots&  \vdots \\
  O&  O&  \cdots&   O&  O&   \cdots&  O&  O&  \cdots&  O&  O \\
  O&  O&  \cdots&   O&  O&    \cdots&  O&  O&  \cdots&  O&  0
\end{bmatrix}$.
\end{center}
Therefore
\begin{equation*}
\begin{aligned}
\lambda_{1}(G)-\lambda_{1}(G_{3})&\ge Y^{T}(D(G)-D(G_{3}))Y\\
&=2\big((n_{q}-h-1)n_{j}y_{j}^{2}-(n_{q}-h-1)\delta y_{q}'y_{j}-(n_{q}-h-1)(h-\delta)y_{q}y_{j}\big)\\
&=2(n_{q}-h-1)y_{j}\big(n_{j}y_{j}-\delta y_{q}'-(h-\delta)y_{q}\big)\\
&>0.
\end{aligned}
\end{equation*}
Then $\lambda_{1}(G_{3})<\lambda_{1}(G)$, which leads a contradiction. Thus $u\in V(B_{i})$ for some $n_{i}=h+1$, where $1\leq i\leq r$. 

In what follows, we shall prove that $G\cong G_{n,(h+1)^{r-1}}^{0,\delta}$. In fact, by the minimality of $\lambda_{1}(G)$ and Claims \ref{claim1.1}-\ref{claim1.3}, we get $G-u\cong K_{c\kappa_{r}^{h}}\vee (K_{n_{1}}\cup \cdots \cup K_{n_{i}-1}\cup \cdots \cup K_{n_{r}})$ for $n_{1}\leq \cdots \leq n_{i} \leq \cdots \leq n_{r}$, where $1\leq i\leq r$. Then $h+1=n_{1}= \cdots = n_{i} \leq \cdots \leq n_{r}$.

\begin{claim}
$n_{a}=h+1$ for $i+1\leq a\leq r-1$.     
\end{claim}

The result follows if $n_{r}=h+1$ or $n_{r-1}=h+1$. Next, consider the case $n_{r-1}>h+1$, so $1\leq i\leq r-2$. Suppose to the contrary that $n_{a}=h+1$. Then there exists some $n_{j}> h+1$ for $i+1\leq j\leq r-1$. Let $V(K_{n_{j}})=\{w_{1},\ldots,w_{n_{j}}\}$ and $G_{4}=G-\{w_{b}w_{c} \mid 1\leq b\leq h+1, h+2\leq c\leq n_{j}\}+\{w_{c}v \mid h+2\leq c\leq n_{j}, v\in V(K_{n_{r}})\}$. Obviously, $G_{4}\in \mathcal{G}_{n,\delta}^{c\kappa_{r}^{h}}$. Let $X$ be the Perron vector of $D(G_{4})$ and $x(v)$ denote the entry of $X$ corresponding to the vertex $v\in V(G_{4})$. By symmetry, set $x(v)=x_{k}$ for any $v\in V(K_{c\kappa_{r}^{h}})$, $x(v)=x_{l}$ for any $v\in V(K_{h+1})$ where $1\leq l\leq i-1$, $x(v)=x_{i}'$ for any $v\in V(K_{\delta})$, $x(v)=x_{i}$ for any $v\in V(K_{h-\delta})$, $x(v)=x_{m}$ for any $v\in V(K_{n_{m}})$ where $i+1\leq m\leq r-1$ and $m\neq j$, $x(v)=x_{r}$ for any $v\in V(K_{n_{j}+n_{r}-h-1})$ and $x(u)=x_{u}$ for the vertex $u$. Then 
$$X=(\underbrace{x_{k},\dots,x_{k}}_{c\kappa_{r}^{h}},\ldots,\underbrace{x_{i}',\dots,x_{i}'}_{\delta},\underbrace{x_{i},\dots,x_{i}}_{h-\delta},\ldots,\underbrace{x_{j},\dots,x_{j}}_{h+1},\ldots,\underbrace{x_{r},\dots,x_{r}}_{n_{j}+n_{r}-h-1},x_{u})^{T}.$$
By $D(G_{4})X=\lambda_{1}(G_{4})X$, we have
\begin{equation*}
\begin{aligned}
\lambda_{1}(G_{4})x_{r}=&c\kappa_{r}^{h}x_{k}+2\sum_{l=1}^{i-1}(h+1)x_{l}+2x_{u}+2\delta x_{i}'+2(h-\delta)x_{i}+2\sum_{m=i+1,m\neq j}^{r-1}n_{m}x_{m}+2(h+1)x_{j}\\
&+(n_{j}+n_{r}-h-2)x_{r},
\end{aligned}
\end{equation*}
\begin{equation*}
\begin{aligned}
\lambda_{1}(G_{4})x_{j}=&c\kappa_{r}^{h}x_{k}+2\sum_{l=1}^{i-1}(h+1)x_{l}+2x_{u}+2\delta x_{i}'+2(h-\delta)x_{i}+2\sum_{m=i+1,m\neq j}^{r-1}n_{m}x_{m}+hx_{j}\\
&+2(n_{j}+n_{r}-h-1)x_{r}.
\end{aligned}
\end{equation*}
Then
\begin{equation*}
\begin{aligned}
&\quad (\lambda_{1}(G_{4})+n_{1}+1)\big(n_{r}x_{r}-(h+1)x_{j}\big)\\
& =(n_{r}-h-1)c\kappa_{r}^{h}x_{k}+\big((n_{r}-2(h+1))(n_{j}+n_{r}-h-1)+n_{r}^{2}\big)x_{r}+(n_{r}-h-1)(h+1)x_{j}\\
& \quad +2(n_{r}-h-1)\Big(\sum_{l=1}^{i-1}(h+1)x_{l}+x_{u}+\delta x_{i}'+(h-\delta)x_{i}+\sum_{m=i+1,m\neq j}^{r-1}n_{m}x_{m}\Big)\\
& >(n_{r}-h-1)c\kappa_{r}^{h}x_{k}+2n_{r}(n_{r}-h-1)x_{r}+(n_{r}-h-1)(h+1)x_{j}\\
& \quad +2(n_{r}-h-1)\Big(\sum_{l=1}^{i-1}(h+1)x_{l}+x_{u}+\delta x_{i}'+(h-\delta)x_{i}+\sum_{m=i+1,m\neq j}^{r}n_{m}x_{m}\Big)\\
& >0.
\end{aligned}
\end{equation*}
Note that $D(G)-D(G_{4})$ is\\
\begin{center}
$\begin{matrix}
c\kappa_{r}^{h}\\
h+1\\
\vdots\\
\delta\\
h-\delta\\
\vdots\\
h+1\\
\vdots\\
n_{j}-h-1\\
n_{r}\\
1
\end{matrix}
\begin{bmatrix}
  O&  O&  \cdots&  O&  O&  \cdots&  O&  \cdots&  O&  O&  O \\
  O&  O&  \cdots&  O&  O&  \cdots&  O&  \cdots&  O&  O&  O \\
  \vdots&  \vdots&  \ddots&  \vdots&  \vdots&  \ddots&  \vdots&  \ddots&  \vdots&  \vdots&  \vdots \\
  O&  O&  \cdots&  O&  O&  \cdots&  O&  \cdots&  O&  O&  O \\
  O&  O&  \cdots&  O&  O&  \cdots&  O&  \cdots&  O&  O&  O \\
  \vdots&  \vdots&  \ddots&  \vdots&  \vdots&  \ddots&  \vdots&  \ddots&  \vdots&  \vdots&  \vdots \\
  O&  O&  \cdots&  O&  O&  \cdots&  O&  \cdots& -J&  O&  O \\
  \vdots&  \vdots&  \ddots&  \vdots&  \vdots&  \ddots&  \vdots&  \ddots&  \vdots&  \vdots&  \vdots \\
  O&  O&  \cdots&  O&  O&  \cdots& -J&  \cdots&  O&  J&  O \\
  O&  O&  \cdots&  O&  O&  \cdots&  O&  \cdots&  J&  O&  O \\
  O&  O&  \cdots&  O&  O&  \cdots&  O&  \cdots&  O&  O&  0
\end{bmatrix}$.
\end{center}
Therefore,
\begin{equation*}
\begin{aligned}
\lambda_{1}(G)-\lambda_{1}(G_{4})&\ge X^{T}(D(G)-D(G_{4}))X\\
&=2\big((n_{j}-h-1)n_{1}x_{1}^{2}-(n_{j}-h-1)(h+1)x_{1}x_{j}\big)\\
&=2(n_{j}-h-1)x_{1}\big(n_{1}x_{1}-(h+1)x_{j}\big)\\
&>0.
\end{aligned}
\end{equation*}
Then $\lambda_{1}(G_{4})<\lambda_{1}(G)$, which leads a contradiction. Thus, $n_{a}=h+1$ for $i+1\leq a\leq r-1$. It suggests that $G\cong G_{n,(h+1)^{r-1}}^{0,\delta}$, as desired.

\begin{case}  
$h<\delta < c\kappa_{r}^{h}+h$.
\end{case}

Similar to Case \ref{case1.1}, we first have the following Claim due to the uncertainty of the minimum degree vertex $u$.

\begin{claim} \label{claim2.1}
$u \notin S$.     
\end{claim}

Otherwise, assume $u\in S$. By the minimality of $\lambda_{1}(G)$, it follows that $G-u\cong K_{c\kappa_{r}^{h}-1}\vee(K_{n_{1}}\cup K_{n_{2}}\cup \cdots \cup K_{n_{r}})$. By using a similar analysis as Claim \ref{claim1.1}, we deduce that $d_{S}(u)=t\leq \delta-2$, $N_{G}(u)\nsubseteq V(S\cup B_{j})$ and $N_{G}(u)\nsubseteq V(B_{j})$ for $1\leq j\leq r$. Then $\sum_{j=1}^{r}p_{j}=\delta-t$. Without loss of generality, suppose that $p_{1}\geq 1$. Depending on the relationship between $n_{1}-p_{1}$ and $\sum_{j=2}^{r}p_{j}$, we divide the subsequent discussion into the following two subcases.

First, consider the case where $n_{1}-p_{1}\geq \sum_{j=2}^{r}p_{j}$. Suppose that $P_{1}=\{w_{1},\ldots,w_{p_{1}}\}$ and $V(K_{n_{1}})\backslash \\ P_{1}=\{w_{p_{1}+1},\ldots,w_{n_{1}}\}$. Let $G_{5}=G-\{uz \mid z\in P_{j}, 2\leq j\leq r\}+\{uw_{a} \mid p_{1}+1\leq a\leq \sum_{j=1}^{r}p_{j}\}+\{w_{p_{1}}v \mid v\in V(K_{n_{j}}), 2\leq j\leq r\}$. Obviously, $G_{5} \in \mathcal{G}_{n,\delta}^{c\kappa_{r}^{h}}$. Let $Y$ be the Perron vector of $D(G_{5})$ and $y(v)$ denote the entry of $Y$ corresponding to the vertex $v\in V(G_{5})$. By symmetry, set $y(v)=y_{d}$ for any $v\in V(K_{t+1})$, $y(v)=y_{k}$ for any $v\in V(K_{c\kappa_{r}^{h}-t-1})$, $y(v)=y_{1}'$ for any $v\in V(K_{\sum_{j=1}^{r}p_{j}-1})$, $y(v)=y_{1}$ for any $v\in V(K_{n_{1}-\sum_{j=1}^{r}p_{j}})$, $y(v)=y_{j}$ for any $v\in V(K_{n_{j}})$ where $2\leq j\leq r$ and $y(u)=y_{u}$ for the vertex $u$. Then 
$$Y=(\underbrace{y_{d},\dots,y_{d}}_{t+1},\underbrace{y_{k},\dots,y_{k}}_{c\kappa_{r}^{h}-t-1},\underbrace{y_{1}',\dots,y_{1}'}_{\sum_{j=1}^{r}p_{j}-1},\underbrace{y_{1},\dots,y_{1}}_{n_{1}-\sum_{j=1}^{r}p_{j}},\ldots,\underbrace{y_{j},\dots,y_{j}}_{n_{j}},\ldots,y_{u})^{T}.$$
By $D(G_{5})Y=\lambda_{1}(G_{5})Y$, we have
$$\lambda_{1}(G_{5})y_{u}=(t+1)y_{d}+2(c\kappa_{r}^{h}-t-1)y_{k}+\Big(\sum_{j=1}^{r}p_{j}-1\Big)y_{1}'+2\Big(n_{1}-\sum_{j=1}^{r}p_{j}\Big)y_{1}+2\sum_{j=2}^{r}n_{j}y_{j},$$
$$\lambda_{1}(G_{5})y_{t}=ty_{d}+(c\kappa_{r}^{h}-t-1)y_{k}+y_{u}+\Big(\sum_{j=1}^{r}p_{j}-1\Big)y_{1}'+\Big(n_{1}-\sum_{j=1}^{r}p_{j}\Big)y_{1}+\sum_{j=2}^{r}n_{j}y_{j},$$
$$\lambda_{1}(G_{5})y_{1}'=(t+1)y_{d}+(c\kappa_{r}^{h}-t-1)y_{k}+y_{u}+\Big(\sum_{j=1}^{r}p_{j}-2\Big)y_{1}'+\Big(n_{1}-\sum_{j=1}^{r}p_{j}\Big)y_{1}+2\sum_{j=2}^{r}n_{j}y_{j},$$
\begin{equation*}
\begin{aligned}
\lambda_{1}(G_{5})y_{j}&=(t+1)y_{d}+(c\kappa_{r}^{h}-t-1)y_{k}+2y_{u}+2\Big(\sum_{j=1}^{r}p_{j}-1\Big)y_{1}'+2\Big(n_{1}-\sum_{j=1}^{r}p_{j}\Big)y_{1}+2\sum_{l=2,l\neq j}^{r}n_{l}y_{l}\\
& \quad +(n_{j}-1)y_{j}.
\end{aligned}
\end{equation*}
Combining the above equalities, we have
\begin{equation} \label{equ4}
\begin{aligned}
 y_{u}-y_{d}=\frac{(c\kappa_{r}^{h}-t-1)y_{k}+(n_{1}-\sum_{j=1}^{r}p_{j})y_{1}+\sum_{j=2}^{r}n_{j}y_{j}}{\lambda_{1}(G_{5})+1}>0,
\end{aligned}
\end{equation}
\begin{equation} \label{equ5}
\begin{aligned}
 y_{1}'-(y_{u}-y_{d})=\frac{ty_{d}+2y_{u}+(\sum_{j=1}^{r}p_{j}-2)y_{1}'+\sum_{j=2}^{r}n_{j}y_{j}}{\lambda_{1}(G_{5})}>0
\end{aligned}
\end{equation} 
and
\begin{equation} \label{equ6}
\begin{aligned}
&\quad (\lambda_{1}(G_{5})+n_{j}+1)(y_{u}-y_{j})\\
& =(c\kappa_{r}^{h}-t-1)y_{k}+(n_{j}-1)y_{u}-\Big(\sum_{j=1}^{r}p_{j}-1\Big)y_{1}'\\
& =(c\kappa_{r}^{h}-t-1)y_{k}+(n_{j}-1)\left(y_{1}'+\frac{(c\kappa_{r}^{h}-t-1)y_{k}+(n_{1}-\sum_{j=1}^{r}p_{j})y_{1}}{\lambda_{1}(G_{5})+1} \right)-\Big(\sum_{j=1}^{r}p_{j}-1\Big)y_{1}'\\
& \geq(c\kappa_{r}^{h}-t-1)y_{k}+(n_{j}-1)y_{1}'-\Big(\sum_{j=1}^{r}p_{j}-1\Big)y_{1}'\\
& =(c\kappa_{r}^{h}-t-1)y_{k}+\Big(n_{j}-\sum_{j=1}^{r}p_{j}\Big)y_{1}'\\
& \geq 0.
\end{aligned}
\end{equation}
Note that $D(G)-D(G_{5})$ is\\
\begin{center}
$\begin{matrix}
1\\
t\\
c\kappa_{r}^{h}-t-1\\
p_{1}-1\\
\sum_{j=2}^{r}p_{j}\\
n_{1}-\sum_{j=1}^{r}p_{j}\\
p_{2}\\
n_{2}-p_{2}\\
\vdots\\
p_{r}\\
n_{r}-p_{r}\\
1
\end{matrix}
\begin{bmatrix}
  0&  O&  O&  O&  O&  O&  J&  J&  \cdots&  J&  J&  0 \\
  O&  O&  O&  O&  O&  O&  O&  O&  \cdots&  O&  O&  O \\
  O&  O&  O&  O&  O&  O&  O&  O&  \cdots&  O&  O&  O \\
  O&  O&  O&  O&  O&  O&  O&  O&  \cdots&  O&  O&  O \\
  O&  O&  O&  O&  O&  O&  O&  O&  \cdots&  O&  O&  J \\
  O&  O&  O&  O&  O&  O&  O&  O&  \cdots&  O&  O&  O \\
  J&  O&  O&  O&  O&  O&  O&  O&  \cdots&  O&  O& -J \\
  J&  O&  O&  O&  O&  O&  O&  O&  \cdots&  O&  O&  O \\
  \vdots&  \vdots&  \vdots&  \vdots&  \vdots&  \vdots&  \vdots&  \vdots&  \ddots&  \vdots&  \vdots&  \vdots \\
  J&  O&  O&  O&  O&  O&  O&  O&  \cdots&  O&  O& -J\\
  J&  O&  O&  O&  O&  O&  O&  O&  \cdots&  O&  O&  O\\
  0&  O&  O&  O&  J&  O& -J&  O&  \cdots& -J&  O&  0
\end{bmatrix}$.
\end{center}
Combining this with (\ref{equ4})-(\ref{equ6}), we get
\begin{equation*}
\begin{aligned}
\lambda_{1}(G)-\lambda_{1}(G_{5})&\ge Y^{T}(D(G)-D(G_{5}))Y\\
&=2\Big(\sum_{j=2}^{r}n_{j}y_{j}y_{d}+\sum_{j=2}^{r}p_{j}y_{1}'y_{u}-\sum_{j=2}^{r}p_{j}y_{j}y_{u}\Big)\\
&\geq2\sum_{j=2}^{r}p_{j}\big(y_{1}'y_{u}-y_{j}(y_{u}-y_{d})\big)~~(\text{since}~n_{j}\geq p_{j})\\
&>~~2\sum_{j=2}^{r}p_{j}(y_{u}-y_{d})(y_{u}-y_{j})(\text{since}~y_{1}'>y_{u}-y_{d})\\
& >0~~(\text{since}~y_{u}>y_{1}~\text{and}~y_{u}\geq y_{j}).
\end{aligned}
\end{equation*}
It follows that $\lambda_{1}(G_{5})<\lambda_{1}(G)$, a contradiction.

Now, consider the remaining case where $n_{1}-p_{1}< \sum_{j=2}^{r}p_{j}$. Let $N_{S}(u)=\{v_{1},\ldots,v_{\delta-\sum_{j=1}^{r}p_{j}}\}$, $S\backslash N_{S}[u]=\{v_{\delta-\sum_{j=1}^{r}p_{j}+1},\ldots,v_{c\kappa_{r}^{h}-1}\}$ and $P_{1}=\{w_{1},\ldots,w_{p_{1}}\}$. Assume that $E_{1}=\{uz \mid z\in P_{j}, 2\leq j\leq r\}$ and $E_{2}=\{uw \mid w\in V(K_{n_{1}})\backslash P_{1}\}+\{uv_{i} \mid \delta-\sum_{j=1}^{r}p_{j}+1\leq i\leq \delta-n_{1}\}+\{w_{p_{1}}v \mid v\in V(K_{n_{j}}), 2\leq j \leq r\}$. Let $G_{6}=G-E_{1}+E_{2}$. Obviously, $G_{6}\in \mathcal{G}_{n,\delta}^{c\kappa_{r}^{h}}$. Let $X$ be the Perron vector of $D(G_{6})$ and $x(v)$ denote the entry of $X$ corresponding to the vertex $v\in V(G_{6})$. By symmetry, set $x(v)=x_{d}$ for any $v\in V(K_{\delta-n_{1}+1})$, $x(v)=x_{k}$ for any $v\in V(K_{c\kappa_{r}^{h}+n_{1}-\delta-1})$, $x(v)=x_{1}$ for any $v\in V(K_{n_{1}-1})$, $x(v)=x_{j}$ for any $v\in V(K_{n_{j}})$ where $2\leq j\leq r$ and $x(u)=x_{u}$ for the vertex $u$. Then
$$X=(\underbrace{x_{d},\dots,x_{d}}_{\delta-n_{1}+1},\underbrace{x_{k},\dots,x_{k}}_{c\kappa_{r}^{h}+n_{1}-\delta-1},\underbrace{x_{1},\dots,x_{1}}_{n_{1}-1},\underbrace{x_{2},\dots,x_{2}}_{n_{2}},\ldots,\underbrace{x_{j},\dots,x_{j}}_{n_{j}},\ldots,x_{u})^{T}.$$
By $D(G_{6})Y=\lambda_{1}(G_{6})Y$, we have
$$\lambda_{1}(G_{6})x_{1}=(\delta-n_{1}+1)x_{d}+(c\kappa_{r}^{h}+n_{1}-\delta-1)x_{k}+x_{u}+(n_{1}-2)x_{1}+2\sum_{j=2}^{r}n_{j}x_{j},$$
$$\lambda_{1}(G_{6})x_{t}=(\delta-n_{1})x_{d}+(c\kappa_{r}^{h}+n_{1}-\delta-1)x_{k}+x_{u}+(n_{1}-1)x_{1}+\sum_{j=2}^{r}n_{j}x_{j},$$
$$\lambda_{1}(G_{6})x_{u}=(\delta-n_{1}+1)x_{d}+2(c\kappa_{r}^{h}+n_{1}-\delta-1)x_{k}+(n_{1}-1)x_{1}+2\sum_{j=2}^{r}n_{j}x_{j},$$
$$\lambda_{1}(G_{6})x_{j}=(\delta-n_{1}+1)x_{d}+(c\kappa_{r}^{h}+n_{1}-\delta-1)x_{k}+2x_{u}+2(n_{1}-1)x_{1}+2\sum_{l=2,l\neq j}^{r}n_{l}x_{l}+(n_{j}-1)x_{j}.$$
Combining the above equalities, we have
\begin{equation} \label{equ7}
\begin{aligned}
x_{1}-x_{d}=\frac{\sum_{j=2}^{r}n_{j}x_{j}}{\lambda_{1}(G_{6})+1}>0,
\end{aligned}
\end{equation}
\begin{equation} \label{equ8}
\begin{aligned}
x_{d}-(x_{u}-x_{d})=\frac{(\delta-n_{1}-1)x_{d}+2x_{u}+(n_{1}-1)x_{1}}{\lambda_{1}(G_{6})}>0
\end{aligned}
\end{equation}
and
\begin{equation} \label{equ9}
\begin{aligned}
x_{u}-x_{1}=\frac{(c\kappa_{h}+n_{1}-\delta-1)x_{k}}{\lambda_{1}(G_{6})+1}>0.
\end{aligned}
\end{equation}
Furthermore,
\begin{equation} \label{equ10}
\begin{aligned}
& \quad (\lambda_{1}(G_{6})+n_{j}+1)(x_{u}-x_{j})\\
&=(c\kappa_{h}+n_{1}-\delta-1)x_{k}+(n_{j}-1)x_{u}-(n_{1}-1)x_{1}\\
&>(c\kappa_{h}+n_{1}-\delta-1)x_{k}+(n_{j}-1)x_{1}-(n_{1}-1)x_{1}~~(\text{by}~(\ref{equ9})) \\
&=(c\kappa_{h}+n_{1}-\delta-1)x_{k}+(n_{j}-n_{1})x_{1}\\
&>0.
\end{aligned}
\end{equation}
Note that $D(G)-D(G_{6})$ is\\
\begin{center}
$\begin{matrix}
1\\
\delta-\sum_{j=1}^{r}p_{j}\\
\sum_{j=1}^{r}p_{j}-n_{1}\\
c\kappa_{r}^{h}+n_{1}-\delta-1\\
p_{1}-1\\
n_{1}-p_{1}\\
p_{2}\\
n_{2}-p_{2}\\
\vdots\\
p_{r}\\
n_{r}-p_{r}\\
1
\end{matrix}
\begin{bmatrix}
  0&  O&  O&  O&  O&  O&  J&  J&  \cdots&  J&  J&  0 \\
  O&  O&  O&  O&  O&  O&  O&  O&  \cdots&  O&  O&  O \\
  O&  O&  O&  O&  O&  O&  O&  O&  \cdots&  O&  O&  J \\
  O&  O&  O&  O&  O&  O&  O&  O&  \cdots&  O&  O&  O \\
  O&  O&  O&  O&  O&  O&  O&  O&  \cdots&  O&  O&  O \\
  O&  O&  O&  O&  O&  O&  O&  O&  \cdots&  O&  O&  J \\
  J&  O&  O&  O&  O&  O&  O&  O&  \cdots&  O&  O& -J \\
  J&  O&  O&  O&  O&  O&  O&  O&  \cdots&  O&  O&  O \\
  \vdots&  \vdots&  \vdots&  \vdots&  \vdots&  \vdots&  \vdots&  \vdots&  \ddots&  \vdots&  \vdots&  \vdots \\
  J&  O&  O&  O&  O&  O&  O&  O&  \cdots&  O&  O& -J\\
  J&  O&  O&  O&  O&  O&  O&  O&  \cdots&  O&  O&  O\\
  0&  O&  J&  O&  O&  J& -J&  O&  \cdots& -J&  O&  0
\end{bmatrix}.$
\end{center}
Combining this with (\ref{equ7})-(\ref{equ10}), we get
\begin{equation*}
\begin{aligned}
\lambda_{1}(G)-\lambda_{1}(G_{6})&\ge X^{T}(D(G)-D(G_{6}))X\\
&=2\Big(\sum_{j=2}^{r}n_{j}x_{j}x_{d}+\Big(\sum_{j=1}^{r}p_{j}-n_{1}\Big)x_{d}x_{u}+(n_{1}-p_{1})x_{1}x_{u}-\sum_{j=2}^{r}p_{j}x_{j}x_{u}\Big)\\
&>2\Big(\sum_{j=2}^{r}n_{j}x_{j}x_{d}+\Big(\sum_{j=1}^{r}p_{j}-n_{1}\Big)x_{d}x_{u}+(n_{1}-p_{1})x_{d}x_{u}-\sum_{j=2}^{r}p_{j}x_{j}x_{u}\Big)\\~~
&\geq 2\sum_{j=2}^{r}p_{i}\big(x_{d}x_{u}-x_{j}(x_{u}-x_{d})\big) ~~(\text{since}~n_{j}\geq p_{j})\\
&>2\sum_{j=2}^{r}p_{j}(x_{u}-x_{d})(x_{u}-x_{j})~~(\text{since}~x_{t}>x_{u}-x_{d})\\
&>0~~(\text{since}~x_{u}>x_{1}~\text{and}~x_{u}>x_{j}).
\end{aligned}
\end{equation*}
Thus, $\lambda_{1}(G_{6})<\lambda_{1}(G)$, which contradicts the minimality of $\lambda_{1}(G)$. This demonstrates that $u \notin S$, proving Claim \ref{claim2.1}.

\vspace{3mm}
Note that $u \notin S$, and so $u \in V(B_{i})$ for some $1\leq i\leq r$ by Claim \ref{claim2.1}. For convenience of notation and calculation, we assume without loss of generality that $u \in V(B_{1})$. Then $d_{B_{1}}(u)\geq 1$ because $n_{1}\geq h+1\geq 2$ and $B_{1}$ is connected.

\begin{claim} \label{claim 2.2}
$t=\delta-h$.     
\end{claim}

If not, $t \neq \delta-h$. This implies that either $d_{B_{1}}(u)<h$ or $d_{B_{1}}(u)>h$. By the minimality of $\lambda_{1}(G)$, we can deduce that $G-u\cong K_{c\kappa_{r}^{h}}\vee (K_{n_{1}-1}\cup K_{n_{2}}\cup \cdots \cup K_{n_{r}})$. If $d_{B_{1}}(u)<h$, let $N_{S}(u)=\{v_{1},\ldots,v_{\delta-p_{1}}\}$, $S\backslash N_{S}(u)=\{v_{\delta-p_{1}+1},\ldots,v_{c\kappa_{r}^{h}}\}$, $P_{1}=\{w_{1},\ldots,w_{p_{1}}\}$ and $V(K_{n_{1}-1})\backslash P_{1}=\{w_{p+1},\ldots,\\w_{n_{1}-1}\}$. Suppose that $G_{7}=G-\{uv_{l} \mid \delta-h+1 \leq l \leq \delta-p_{1}\}+\{uw_{m} \mid p+1 \leq m \leq h\}$. Obviously, $G_{7} \in \mathcal{G}_{n,\delta}^{c\kappa_{r}^{h}}$. Let $Y$ be the Perron vector of $D(G_{7})$ and $y(v)$ denote the entry of $Y$ corresponding to the vertex $v\in V(G_{7})$. By symmetry, set $y(v)=y_{d}$ for any $v\in V(K_{\delta-h})$, $y(v)=y_{k}$ for any $v\in V(K_{c\kappa_{r}^{h}+h-\delta})$, $y(v)=y_{1}'$ for any $v\in V(K_{h})$, $y(v)=y_{1}$ for any $v\in V(K_{n_{1}-h-1})$, $y(v)=y_{i}$ for any $v\in V(K_{n_{i}})$ where $2\leq i\leq r$ and $y(u)=y_{u}$ for the vertex $u$. Then 
$$Y=(\underbrace{y_{d},\dots,y_{d}}_{\delta-h},\underbrace{y_{k},\dots,y_{k}}_{c\kappa_{r}^{h}+h-\delta},\underbrace{y_{1}',\dots,y_{1}'}_{h},\underbrace{y_{1},\dots,y_{1}}_{n_{1}-h-1},\ldots,\underbrace{y_{i},\dots,y_{i}}_{n_{i}},\ldots,y_{u})^{T}.$$
By $D(G_{7})Y=\lambda_{1}(G_{7})Y$, we have
$$\lambda_{1}(G_{7})y_{1}'=(\delta-h)y_{d}+(c\kappa_{r}^{h}+h-\delta)y_{k}+y_{u}+(h-1)y_{1}'+(n_{1}-h-1)y_{1}+2\sum_{i=2}^{r}n_{i}y_{i},$$
$$\lambda_{1}(G_{7})y_{k}=(\delta-h)y_{d}+(c\kappa_{r}^{h}+h-\delta-1)y_{k}+2y_{u}+hy_{1}'+(n_{1}-h-1)y_{1}+\sum_{i=2}^{r}n_{i}y_{i},$$
$$\lambda_{1}(G_{7})y_{i}=(\delta-h)y_{d}+(c\kappa_{r}^{h}+h-\delta)y_{k}+2y_{u}+2hy_{1}'+2(n_{1}-h-1)y_{1}+\sum_{j=2,j\neq i}^{r}n_{j}y_{j}+(n_{i}-1)y_{i},$$
$$\lambda_{1}(G_{7})y_{u}=(\delta-h)y_{d}+2(c\kappa_{r}^{h}+h-\delta)y_{k}+hy_{1}'+2(n_{1}-h-1)y_{1}+2\sum_{i=2}^{r}n_{i}y_{i}.$$
From which we get
\begin{equation*}
\begin{aligned}
& \quad (\lambda_{1}(G_{7})+3)\Big(\sum_{i=2}^{r}n_{i}y_{i}-y_{u}\Big)\\
&=\Big(\sum_{i=2}^{r}n_{i}-1\Big)(\delta-h)y_{d}+\Big(\sum_{i=2}^{r}n_{i}-2\Big)(c\kappa_{r}^{h}+h-\delta)y_{k}+\Big(2\sum_{i=2}^{r}n_{i}-3\Big)y_{u}+\Big(2\sum_{i=2}^{r}n_{i}-1\Big)hy_{1}'\\
& \quad +2\Big(\sum_{i=2}^{r}n_{i}-1\Big)(n_{1}-h-1)y_{1}+\sum_{i=2}^{r}\sum_{j=2,j\neq i}^{r}n_{i}n_{j}y_{j}+\sum_{i=2}^{r}n_{i}^{2}y_{i}\\
&>0.
\end{aligned}
\end{equation*}
Then
$$y_{1}'-y_{k}=\frac{\sum_{i=2}^{r}n_{i}y_{i}-y_{u}}{\lambda_{1}(G_{7})+1}>0.$$
Note that $D(G)-D(G_{7})$ is\\
\begin{center}
$\begin{matrix}
\delta-h\\
h-p_{1}\\
c\kappa_{r}^{h}+p_{1}-\delta\\
p_{1}\\
h-p_{1}\\
n_{1}-h-1\\
n_{2}\\
\vdots\\
n_{r}\\
1
\end{matrix}
\begin{bmatrix}
  O&  O&  O&  O&  O&  O&  O&  \cdots&  O&  O \\
  O&  O&  O&  O&  O&  O&  O&  \cdots&  O& -J \\
  O&  O&  O&  O&  O&  O&  O&  \cdots&  O&  O \\
  O&  O&  O&  O&  O&  O&  O&  \cdots&  O&  O \\
  O&  O&  O&  O&  O&  O&  O&  \cdots&  O&  J \\
  O&  O&  O&  O&  O&  O&  O&  \cdots&  O&  O\\
  O&  O&  O&  O&  O&  O&  O&  \cdots&  O&  O\\
  \vdots&  \vdots&  \vdots&  \vdots&  \vdots&  \vdots&  \vdots&  \ddots&  \vdots&  \vdots \\
  O&  O&  O&  O&  O&  O&  O&  \cdots&  O&  O\\
  O& -J&  O&  O&  J&  O&  O&  \cdots&  O&  0
\end{bmatrix}$.
\end{center}
Therefore
\begin{equation*}
\begin{aligned}
\lambda_{1}(G)-\lambda_{1}(G_{7})&\ge Y^{T} (D(G)-D(G_{7}))Y\\
&=2\big((h-p_{1})y_{1}'y_{u}-(h-p_{1})y_{k}y_{u}\big)\\
&=2(h-p_{1})y_{u}(y_{1}'-y_{k})\\
&>0~~(\text{since}~ y_{1}'>y_{k}).
\end{aligned}
\end{equation*}
Hence $\lambda_{1}(G_{7})<\lambda_{1}(G)$, a contradiction.

If $d_{B_{1}}(u)>h$, let $N_{S}(u)=\{v_{1},\ldots,v_{\delta-p_{1}}\}$, $S\backslash N_{S}(u)=\{v_{\delta-p_{1}+1},\ldots,v_{c\kappa_{r}^{h}}\}$, $P_{1}=\{w_{1},\ldots,w_{p_{1}}\}$ and $V(K_{n_{1}-1})\backslash P_{1}=\{w_{p+1},\ldots,w_{n_{1}-1}\}$. Let $E_{3}=\{uw_{a} \mid h+1 \leq a \leq p_{1}\}+\{w_{a}w_{b} \mid 1 \leq a \leq h, h+1 \leq b\leq n_{1}-1\}$ and $E_{4}=\{uv_{c} \mid \delta-p_{1}+1 \leq c \leq\delta-h\}+\{w_{b}z \mid h+1 \leq b \leq n_{1}-1, z\in V(K_{n_{2}})\}$. Then $G_{8}=G-E_{3}+E_{4}$. Obviously, $G_{8}\in \mathcal{G}_{n,\delta}^{c\kappa_{r}^{h}}$. Let $X$ be the Perron vector of $D(G_{8})$ and $x(v)$ denote the entry of $X$ corresponding to the vertex $v\in V(G_{8})$. By symmetry, set $x(v)=x_{d}$ for any $v\in V(K_{\delta-h})$, $x(v)=x_{k}$ for any $v\in V(K_{c\kappa_{r}^{h}+h-\delta})$, $x(v)=x_{1}$ for any $v\in V(K_{h})$, $x(v)=x_{2}$ for any $v\in V(K_{n_{1}+n_{2}-h-1})$,$x(v)=x_{i}$ for any $v\in V(K_{n_{i}})$ where $3\leq i\leq r$ and $x(u)=x_{u}$ for the vertex $u$. Then
$$Y=(\underbrace{x_{d},\dots,x_{d}}_{\delta-h},\underbrace{x_{k},\dots,x_{k}}_{c\kappa_{r}^{h}+h-\delta},\underbrace{x_{1},\dots,x_{1}}_{h},\underbrace{x_{2},\dots,x_{2}}_{n_{1}+n_{2}-h-1},\ldots,\underbrace{x_{i},\dots,x_{i}}_{n_{i}},\ldots,x_{u})^{T}.$$
By $D(G_{8})X=\lambda_{1}(G_{8})X$, we have
$$\lambda_{1}(G_{8})x_{d}=(\delta-h-1)x_{d}+(c\kappa_{r}^{h}+h-\delta)x_{k}+x_{u}+hx_{1}+(n_{1}+n_{2}-h-1)x_{2}+\sum_{i=3}^{r}n_{i}x_{i},$$
$$\lambda_{1}(G_{8})x_{2}=(\delta-h)x_{d}+(c\kappa_{r}^{h}+h-\delta)x_{k}+2x_{u}+2hx_{1}+(n_{1}+n_{2}-h-2)x_{2}+2\sum_{i=3}^{r}n_{i}x_{i},$$
$$\lambda_{1}(G_{8})x_{1}=(\delta-h)x_{d}+(c\kappa_{r}^{h}+h-\delta)x_{k}+x_{u}+(h-1)x_{1}+2(n_{1}+n_{2}-h-1)x_{2}+2\sum_{i=3}^{r}n_{i}x_{i},$$
$$\lambda_{1}(G_{8})x_{u}=(\delta-h)x_{d}+2(c\kappa_{r}^{h}+h-\delta)x_{k}+hx_{1}+2(n_{1}+n_{2}-h-1)x_{2}+2\sum_{i=3}^{r}n_{i}x_{i},$$
from which we get
\begin{equation} \label {equ11}
\begin{aligned}
&(\lambda_{1}(G_{8})+1)(n_{2}x_{d}-x_{u})\\
=&(n_{2}-1)(\delta-h)x_{d}+(n_{2}-2)(c\kappa_{r}^{h}+h-\delta)x_{k}+(n_{2}-1)x_{u}+(n_{2}-2)hx_{1}+(n_{2}-2)\sum_{i=3}^{r}n_{i}x_{i}\\
& +(n_{2}-2)(n_{1}+n_{2}-h-1)x_{2}\\
>&0,
\end{aligned}
\end{equation}
\begin{equation} \label {equ12}
\begin{aligned}
2x_{2}-x_{u}=\frac{(\delta-h)x_{d}+3x_{u}+3hx_{1}+2\sum_{i=3}^{r}n_{i}x_{i}}{\lambda_{1}(G_{8})+1}>0
\end{aligned}
\end{equation}
and
\begin{equation} \label{equ13}
\begin{aligned}
&(\lambda_{1}(G_{8})+n_{2}+1)\big((2n_{2}-h)x_{2}-hx_{1}\big)\\
=&2(n_{2}-h)(\delta-h)x_{d}+2(n_{2}-h)(c\kappa_{r}^{h}+h-\delta)x_{k}+4(n_{2}-h)x_{u}+3(n_{2}-h)hx_{1}+4(n_{2}-h)\sum_{i=3}^{r}n_{i}x_{i}\\
& +\big((2n_{2}-3h)(n_{1}+n_{2}-h-1)+n_{2}(2n_{2}-h)\big)x_{2}\\
\geq&2(n_{2}-h)(\delta-h)x_{d}+2(n_{2}-h)(c\kappa_{r}^{h}+h-\delta)x_{k}+4(n_{2}-h)x_{u}+3(n_{2}-h)hx_{1}+4(n_{2}-h)\sum_{i=3}^{r}n_{i}x_{i}\\
& +4n_{2}(n_{2}-h)x_{2}\\
>&0.
\end{aligned}
\end{equation}
By (\ref{equ11})-(\ref{equ13}), we have
\begin{equation*}
\begin{aligned}
&(\lambda_{1}(G_{8})+n_{2}+1)\big((p_{1}-h)x_{u}(x_{d}-x_{2})+(n_{1}-h-1)x_{2}(n_{2}x_{2}-hx_{1})\big)\\
=&(p_{1}-h)x_{u}\Big(n_{2}x_{d}-x_{u}-hx_{1}-\sum_{i=3}^{r}n_{i}x_{i}-n_{2}x_{2}\Big)+(n_{1}-h-1)x_{2}\Big((n_{2}-h)(\delta-h)x_{d}+(2n_{2}-h)x_{u}\\
+&(n_{2}-h)(c\kappa_{r}^{h}+h-\delta)x_{k}+(n_{2}-h)hx_{1}+\big((n_{2}-2h)(n_{1}+n_{2}-h-1)+n_{2}^{2}\big)x_{2}+2(n_{2}-h)\sum_{i=3}^{r}n_{i}x_{i}\Big)\\
>&(p_{1}-h)x_{u}(-hx_{1}-\sum_{i=3}^{r}n_{i}x_{i}-n_{2}x_{2})+(n_{1}-h-1)x_{2}\Big((n_{2}-h)(\delta-h)x_{d}+(n_{2}-h)(c\kappa_{r}^{h}+h-\delta)x_{k}\\
+&(2n_{2}-h)x_{u}+(n_{2}-h)hx_{1}+\big((n_{2}-2h)(n_{1}+n_{2}-h-1)+n_{2}^{2}\big)x_{2}+2(n_{2}-h)\sum_{i=3}^{r}n_{i}x_{i}\Big)\\
\geq&(p_{1}-h)x_{u}(-hx_{1}-\sum_{i=3}^{r}n_{i}x_{i}-n_{2}x_{2})+(n_{1}-h-1)x_{2}\Big((n_{2}-h)(\delta-h)x_{d}+(n_{2}-h)(c\kappa_{r}^{h}+h-\delta)x_{k}\\
+&(2n_{2}-h)x_{u}+(n_{2}-h)hx_{1}+2n_{2}(n_{2}-h)x_{2}+2(n_{2}-h)\sum_{i=3}^{r}n_{i}x_{i}\Big)\\
\geq&(p_{1}-h)\Big(-hx_{1}x_{u}-\sum_{i=3}^{r}n_{i}x_{i}x_{u}-n_{2}x_{2}x_{u}+(n_{2}-h)(\delta-h)x_{2}x_{d}+(n_{2}-h)(c\kappa_{r}^{h}+h-\delta)x_{2}x_{k}\\
+&(2n_{2}-h)x_{2}x_{u}+(n_{2}-h)hx_{1}x_{2}+2n_{2}(n_{2}-h)x_{2}^{2}+2(n_{2}-h)\sum_{i=3}^{r}n_{i}x_{i}x_{2}\Big)\\
>&(p_{1}-h)\Big(-hx_{1}x_{u}-\sum_{i=3}^{r}n_{i}x_{i}x_{u}-n_{2}x_{2}x_{u}+(2n_{2}-h)x_{2}x_{u}+2n_{2}(n_{2}-h)x_{2}^{2}+2(n_{2}-h)\sum_{i=3}^{r}n_{i}x_{i}x_{2}\Big)\\
=&(p_{1}-h)\Big(\sum_{i=3}^{r}n_{i}x_{i}\big(2(n_{2}-h)x_{2}-x_{u}\big)+n_{2}x_{2}\big(2(n_{2}-h)x_{2}-x_{u}\big)+x_{u}\big((2n_{2}-h)x_{2}-hx_{1}\big)\Big)\\
\geq&(p_{1}-h)\Big(\sum_{i=3}^{r}n_{i}x_{i}(2x_{2}-x_{u})+n_{2}x_{2}(2x_{2}-x_{u})+x_{u}\big((2n_{2}-h)x_{2}-hx_{1}\big)\Big)\\
>&0.
\end{aligned}
\end{equation*}

Note that $D(G)-D(G_{8})$ is\\
\begin{center}
$\begin{matrix}
\delta-p_{1}\\
p_{1}-h\\
c\kappa_{r}^{h}+h-\delta\\
h\\
p_{1}-h\\
n_{1}-p_{1}-1\\
n_{2}\\
n_{3}\\
\vdots\\
n_{r}\\
1
\end{matrix}
\begin{bmatrix}
  O&  O&  O&  O&  O&  O&  O&  O&  \cdots&  O&  O \\
  O&  O&  O&  O&  O&  O&  O&  O&  \cdots&  O&  J \\
  O&  O&  O&  O&  O&  O&  O&  O&  \cdots&  O&  O \\
  O&  O&  O&  O&  O&  O&  O&  O&  \cdots&  O&  O \\
  O&  O&  O& -J&  O&  O&  J&  O&  \cdots&  O& -J \\
  O&  O&  O& -J&  O&  O&  J&  O&  \cdots&  O&  O \\
  O&  O&  O&  O&  J&  J&  O&  O&  \cdots&  O&  O \\
  O&  O&  O&  O&  O&  O&  O&  O&  \cdots&  O&  O \\
  \vdots&  \vdots&  \vdots&  \vdots&  \vdots&  \vdots&  \vdots&  \vdots&  \ddots&  \vdots&  \vdots \\
  O&  O&  O&  O&  O&  O&  O&  O&  \cdots&  O&  O \\
  O&  J&  O&  O& -J&  O&  O&  O&  \cdots&  O&  0
\end{bmatrix}$.
\end{center}
Combining the above inequalities, we obtain
\begin{equation*}
\begin{aligned}
\lambda_{1}(G)-\lambda_{1}(G_{8})&\ge X^{T} (D(G)-D(G_{8}))X\\
&=2\big((p_{1}-h)x_{d}x_{u}+(n_{1}-h-1)n_{2}x_{2}^{2}-(p_{1}-h)x_{2}x_{u}-(n_{1}-h-1)hx_{1}x_{2}\big)\\
&=2\big((p_{1}-h)x_{u}(x_{d}-x_{2})+(n_{1}-h-1)(n_{2}x_{2}-hx_{1})\big)\\
&>0.
\end{aligned}
\end{equation*}
Therefore $\lambda_{1}(G_{8})<\lambda_{1}(G)$, a contradiction. Then we conclude that $t=\delta-h$.

\begin{claim} \label{claim2.3}
$u\in V(B_{i})$ for some $n_{i}=h+1$, where $1\leq i\leq r$.     
\end{claim}
If $n_{r}=h+1$, the result is trivial. Next, we consider the case of $n_{r}>h+1$. We use proof by contradiction, assume that there exists $B_{q}$ such that $u\in V(B_{q})$ for $n_{q}>h+1$, where $1\leq q\leq r$. By the minimality of $\lambda_{1}(G)$, we can infer that $G-u\cong K_{c\kappa_{r}^{h}}\vee (K_{n_{1}}\cup \cdots \cup K_{n_{q}-1}\cup \cdots \cup K_{n_{r}})$. Suppose that $P_{q}=\{w_{1},\ldots, w_{h}\}$ and $V(K_{n_{q}-1})\backslash P_{q}=\{w_{h+1},\ldots, w_{n_{q}-1}\}$. Let $G_{9}=G-\{w_{a}w_{b} \mid 1\leq a\leq h, h+1\leq b\leq n_{q}-1\}+\{w_{b}z \mid h+1\leq b\leq n_{q}-1, z\in V(K_{n_{j}})\}$, where $j\neq q$. Obviously, $G_{9}\in \mathcal{G}_{n,\delta}^{c\kappa_{r}^{h}}$. Let $Y$ be the Perron vector of $D(G_{9})$ and $y(v)$ denote the entry of $Y$ corresponding to the vertex $v\in V(G_{9})$. By symmetry, set $y(v)=y_{d}$ for any $v\in V(K_{\delta-h})$, $y(v)=y_{k}$ for any $v\in V(K_{c\kappa_{r}^{h}+h-\delta})$, $y(v)=y_{q}$ for any $v\in V(K_{h})$, $y(v)=y_{j}$ for any $v\in V(K_{n_{q}+n_{j}-h-1})$, $y(v)=y_{l}$ for any $v\in V(K_{n_{l}})$ where $1 \leq l \leq r$ and $l\neq q,j$ and $y(u)=y_{u}$ for the vertex $u$. Then
$$Y=(\underbrace{y_{d},\dots,y_{d}}_{\delta-h},\underbrace{y_{k},\dots,y_{k}}_{c\kappa_{r}^{h}},\ldots,\underbrace{y_{i},\dots,y_{i}}_{h},\ldots,\underbrace{y_{j},\dots,y_{j}}_{n_{i}+n_{j}-h-1},\ldots,y_{u})^{T}.$$
By $D(G_{9})Y=\lambda_{1}(G_{9})Y$, we have
$$\lambda_{1}(G_{9})y_{j}=(\delta-h)y_{d}+(c\kappa_{r}^{h}+h-\delta)y_{k}+2y_{u}+2hy_{q}+(n_{q}+n_{j}-h-2)y_{j}+2\sum_{l=1,l\neq q,j}^{r}n_{l}y_{l},$$
$$\lambda_{1}(G_{9})y_{q}=(\delta-h)y_{d}+(c\kappa_{r}^{h}+h-\delta)y_{k}+y_{u}+(h-1)y_{q}+2(n_{q}+n_{j}-h-1)y_{j}+2\sum_{l=1,l\neq q,j}^{r}n_{l}y_{l}.$$
Then
\begin{equation*}
\begin{aligned}
& \quad (\lambda_{1}(G_{9})+n_{j}+1)(n_{j}y_{j}-hy_{q})\\
& =(n_{j}-h)(\delta-h)y_{d}+(n_{j}-h)(c\kappa_{r}^{h}+h-\delta)y_{k}+(2n_{j}-h)y_{u}+(n_{j}-h)hy_{q}+2(n_{j}-h)\sum_{l=1,l\neq q,j}^{r}n_{l}y_{l}\\
& \quad +\big((n_{j}-2h)(n_{q}+n_{j}-h-1)+n_{j}^{2}\big)y_{j}\\
& \geq (n_{j}-h)(\delta-h)y_{d}+(n_{j}-h)(c\kappa_{r}^{h}+h-\delta)y_{k}+(2n_{j}-h)y_{u}+(n_{j}-h)hy_{q}+2(n_{j}-h)\sum_{l=1,l\neq q,j}^{r}n_{l}y_{l}\\
& \quad +\big((n_{j}-2h)n_{j}+n_{j}^{2}\big)y_{j}\\
& =(n_{j}-h)(\delta-h)y_{d}+(n_{j}-h)(c\kappa_{r}^{h}+h-\delta)y_{k}+(2n_{j}-h)y_{u}+(n_{j}-h)hy_{q}+2(n_{j}-h)\sum_{l=1,l\neq q,j}^{r}n_{l}y_{l}\\
& \quad +2n_{j}(n_{j}-h)y_{j}\\
& >0.
\end{aligned}
\end{equation*}
Note that $D(G)-D(G_{9})$ is\\
\begin{center}
$\begin{matrix}
\delta-h\\
c\kappa_{r}^{h}+h-\delta\\
n_{1}\\
\vdots\\
h\\
\vdots\\
n_{q}-h-1\\
n_{j}\\
\vdots\\
n_{r}\\
1
\end{matrix}
\begin{bmatrix}
  O&  O&  O&  \cdots&  O&  \cdots&  O&  O&  \cdots&  O&  O \\
  O&  O&  O&  \cdots&  O&  \cdots&  O&  O&  \cdots&  O&  O \\
  O&  O&  O&  \cdots&  O&  \cdots&  O&  O&  \cdots&  O&  O \\
  \vdots&  \vdots&  \vdots&  \ddots&  \vdots&  \ddots&  \vdots&  \vdots&  \ddots&  \vdots&  \vdots \\
  O&  O&  O&  \cdots&  O&  \cdots& -J&  O&  \cdots&  O&  O \\
  \vdots&  \vdots&  \vdots&  \ddots&  \vdots&  \ddots&  \vdots&  \vdots&  \ddots&  \vdots&  \vdots \\
  O&  O&  O&  \cdots& -J&  \cdots&   O&  J&  \cdots&  O&  O \\
  O&  O&  O&  \cdots&  O&  \cdots&   J&  O&  \cdots&  O&  O \\
  \vdots&  \vdots&  \vdots&  \ddots&  \vdots&  \ddots&  \vdots&  \vdots&  \ddots&  \vdots&  \vdots \\
  O&  O&  O&  \cdots&  O&  \cdots&  O&  O&  \cdots&  O&  O \\
  O&  O&  O&  \cdots&  O&  \cdots&  O&  O&  \cdots&  O&  0
\end{bmatrix}$.
\end{center}
Therefore,
\begin{align*}
\lambda_{1} (D(G))-\lambda_{1} (D(G_{9}))&\ge Y^{T} (D(G)-D(G_{9}))Y\\
&=2\big((n_{i}-h-1)n_{j}y_{j}^{2}-(n_{i}-h-1)hy_{i}y_{j}\big)\\
&=2(n_{i}-h-1)y_{j}(n_{j}y_{j}-hy_{i})\\
&>0.
\end{align*} 
Then $\lambda_{1}(G_{9})<\lambda_{1}(G)$, which leads a contradiction. Thus $u\in V(B_{i})$ for some $n_{i}=h+1$, where $1\leq i\leq r$.


By the minimality of $\lambda_{1}(G)$ and Claims \ref{claim2.1}-\ref{claim2.3}, we get $G-u\cong K_{c\kappa_{r}^{h}}\vee (K_{n_{1}}\cup \cdots \cup K_{n_{i}-1}\cup \cdots \cup K_{n_{r}})$ for $n_{1}\leq \cdots \leq n_{i} \leq \cdots \leq n_{r}$, where $1\leq i\leq r$. Then $h+1=n_{1}= \cdots = n_{i} \leq \cdots \leq n_{r}$.
\begin{claim}
$n_{a}=h+1$ for $i+1\leq a\leq r-1$.     
\end{claim}

The result follows if $n_{r}=h+1$ or $n_{r-1}=h+1$. Next, consider the case $n_{r-1}>h+1$, so $1\leq i\leq r-2$. Suppose to the contrary that $n_{a}=h+1$. Then there exists some $n_{j}\geq h+2$ for $i+1\leq j\leq r-1$. Let $V(K_{n_{j}})=\{w_{1},\ldots,w_{n_{j}}\}$ and $G_{10}=G-\{w_{b}w_{c} \mid 1\leq b\leq h+1, h+2\leq c\leq n_{j}\}+\{w_{c}v \mid h+2\leq c\leq n_{j}, v\in V(K_{n_{r}})\}$. Obviously, $G_{10}\in \mathcal{G}_{n,\delta}^{c\kappa_{r}^{h}}$. Let $X$ be the Perron vector of $D(G_{10})$ and $x(v)$ denote the entry of $X$ corresponding to the vertex $v\in V(G_{10})$. By symmetry, set $x(v)=x_{d}$ for any $v\in V(K_{\delta-h})$, $x(v)=x_{k}$ for any $v\in V(K_{c\kappa_{r}^{h}+h-\delta})$, $x(v)=x_{l}$ for any $v\in V(K_{h+1})$ where $1\leq l\leq i-1$, $x(v)=x_{i}$ for any $v\in V(K_{h})$, $x(v)=x_{m}$ for any $v\in V(K_{n_{m}})$ where $i+1\leq m\leq r-1$ and $m\neq j$, $x(v)=x_{r}$ for any $v\in V(K_{n_{j}+n_{r}-h-1})$ and $x(u)=x_{u}$ for the vertex $u$. Then 
$$X=(\underbrace{x_{d},\dots,x_{d}}_{\delta-h},\underbrace{x_{k},\dots,x_{k}}_{c\kappa_{r}^{h}+h-\delta},\ldots,\underbrace{x_{i},\dots,x_{i}}_{h},\ldots,\underbrace{x_{j},\dots,x_{j}}_{h+1},\ldots,\underbrace{x_{r},\dots,x_{r}}_{n_{j}+n_{r}-h-1},x_{u})^{T}.$$
By $D(G_{10})X=\lambda_{1}(G_{10})X$, we have
\begin{equation*}
\begin{aligned}
\lambda_{1}(G_{10})x_{1}=&(\delta-h)x_{d}+(c\kappa_{r}^{h}+h-\delta)x_{k}+2\sum_{l=1}^{i-1}(h+1)x_{l}+2x_{u}+2hx_{i}+2\sum_{m=i+1,m\neq j}^{r-1}n_{m}x_{m}\\
&+2(h+1)x_{j}+(n_{j}+n_{r}-h-2)x_{r},
\end{aligned}
\end{equation*}
\begin{equation*}
\begin{aligned}
\lambda_{1}(G_{10})x_{j}=&(\delta-h)x_{d}+(c\kappa_{r}^{h}+h-\delta)x_{k}+2\sum_{l=1}^{i-1}(h+1)x_{l}+2x_{u}+2hx_{i}+2\sum_{m=i+1,m\neq j}^{r-1}n_{m}x_{m}\\
&+hx_{j}+(n_{j}+n_{r}-h-1)x_{r}.
\end{aligned}
\end{equation*}
Then
\begin{equation*}
\begin{aligned}
&\quad (\lambda_{1}(G_{10})+n_{r}+1)\big(n_{r}x_{r}-(h+1)x_{j}\big)\\
& =(n_{r}-h-1)(\delta-h)x_{d}+(n_{r}-h-1)(c\kappa_{r}^{h}+h-\delta)x_{k}+\big((n_{r}-2(h+1))(n_{j}+n_{r}-h-1)+n_{r}^{2}\big)x_{r}\\
& \quad +(n_{r}-h-1)(h+1)x_{j}+2(n_{r}-h-1)\Big(\sum_{l=1}^{i-1}(h+1)x_{l}+x_{u}+hx_{i}+\sum_{m=i+1,m\neq j}^{r-1}n_{m}x_{m}\Big)\\
& >(n_{r}-h-1)(\delta-h)x_{d}+(n_{r}-h-1)(c\kappa_{r}^{h}+h-\delta)x_{k}+\big(n_{r}-2(h+1)+n_{r}^{2}\big)x_{r}\\
& \quad +(n_{r}-h-1)(h+1)x_{j}+2(n_{r}-h-1)\Big(\sum_{l=1}^{i-1}(h+1)x_{l}+x_{u}+hx_{i}+\sum_{m=i+1,m\neq j}^{r-1}n_{m}x_{m}\Big)\\
& =(n_{r}-h-1)(\delta-h)x_{d}+(n_{r}-h-1)(c\kappa_{r}^{h}+h-\delta)x_{k}+2n_{1}(n_{r}-h-1)x_{r}+(n_{r}-h-1)(h+1)x_{j}\\
& \quad +2(n_{r}-h-1)\Big(\sum_{l=1}^{i-1}(h+1)x_{l}+x_{u}+hx_{i}+\sum_{m=i+1,m\neq j}^{r-1}n_{m}x_{m}\Big)\\
& >0.
\end{aligned}
\end{equation*}
Note that $D(G)-D(G_{10})$ is\\
\begin{center}
$\begin{matrix}
\delta-h\\
c\kappa_{r}^{h}+h-\delta\\
h+1\\
\vdots\\
h\\
\vdots\\
h+1\\
\vdots\\
n_{j}+n_{r}-h-1\\
n_{r}\\
1
\end{matrix}
\begin{bmatrix}
  O&  O&  O&  \cdots&  O&  \cdots&  O&  \cdots&  O&  O&  O \\
  O&  O&  O&  \cdots&  O&  \cdots&  O&  \cdots&  O&  O&  O \\
  O&  O&  O&  \cdots&  O&  \cdots&  O&  \cdots&  O&  O&  O \\
  \vdots&  \vdots&  \vdots&  \ddots&  \vdots&  \ddots&  \vdots&  \ddots&  \vdots&  \vdots&  \vdots \\
  O&  O&  O&  \cdots&  O&  \cdots&  O&  \cdots&  O&  O&  O \\
  \vdots&  \vdots&  \vdots&  \ddots&  \vdots&  \ddots&  \vdots&  \ddots&  \vdots&  \vdots&  \vdots \\
  O&  O&  O&  \cdots&  O&  \cdots&  O&  \cdots& -J&  O&  O \\
  \vdots&  \vdots&  \vdots&  \ddots&  \vdots&  \ddots&  \vdots&  \ddots&  \vdots&  \vdots&  \vdots \\
  O&  O&  O&  \cdots&  O&  \cdots& -J&  \cdots&  O&  J&  O \\
  O&  O&  O&  \cdots&  O&  \cdots&  O&  \cdots&  J&  O&  O \\
  O&  O&  O&  \cdots&  O&  \cdots&  O&  \cdots&  O&  O&  0
\end{bmatrix}$.
\end{center}
Therefore,
\begin{equation*}
\begin{aligned}
\lambda_{1}(G)-\lambda_{1}(G_{10})&\ge X^{T}(D(G)-D(G_{10}))X\\
&=2\big((n_{j}-h-1)n_{r}x_{r}^{2}-(n_{j}-h-1)(h+1)x_{r}x_{j}\big)\\
&=2(n_{j}-h-1)x_{r}\big(n_{r}x_{r}-(h+1)x_{j}\big)\\
&>0.
\end{aligned}
\end{equation*}
It follows that $\lambda_{1}(G_{10})<\lambda_{1}(G)$, a contradiction. Thus, $n_{a}=h+1$ for $i+1\leq a\leq r-1$. Then we deduce that $G\cong G_{n,(h+1)^{r-1}}^{\delta-h,h}$, as desired.

\begin{case}  
$\delta\geq c\kappa_{r}^{h}+h$.
\end{case}
Recall that $|V(B_{i})|=n_{i}$ for $1\leq i\leq r$ and $n_{1}\leq \cdots \leq n_{i} \leq \cdots \leq n_{r}$. It is straight forward to derive that
$$\lambda_{1}(K_{c\kappa_{r}^{h}}\cup (K_{n_{1}}\cup K_{n_{2}}\cup \cdots \cup K_{n_{r}}))\leq \lambda_{1}(G),$$
with equality if and only if $G\cong K_{c\kappa_{r}^{h}}\cup (K_{n_{1}}\cup K_{n_{2}}\cup \cdots \cup K_{n_{r}}).$ Obviously, $n_{1}\geq \delta-c\kappa_{r}^{h}+1$ since the minimum degree of $G$ is $\delta$. Combining this with Lemma \ref{lem1.2}, we have
$$\lambda_{1}(K_{c\kappa_{r}^{h}}\vee(K_{n-c\kappa_{r}^{h}-(r-1)(\delta-c\kappa_{r}^{h}+1)}\cup (r-1)K_{\delta-c\kappa_{r}^{h}+1})) \leq \lambda_{1}(K_{c\kappa_{r}^{h}}\cup (K_{n_{1}}\cup K_{n_{2}}\cup \cdots \cup K_{n_{r}})),$$
with equality if and only if $(n_{1},\ldots,n_{r-1},n_{r})=(\delta-c\kappa_{r}^{h}+1,\ldots,\delta-c\kappa_{r}^{h}+1,n-c\kappa_{r}^{h}-(r-1)(\delta-c\kappa_{r}^{h}+1))$. By the minimality of $\lambda_{1}(G)$, we conclude that $G\cong (K_{c\kappa_{r}^{h}}\vee(K_{n-c\kappa_{r}^{h}-(r-1)(\delta-c\kappa_{r}^{h}+1)}\cup (r-1)K_{\delta-c\kappa_{r}^{h}+1})$.

This completes the proof.
\begin{flushright}
$\square$
\end{flushright}

\end{document}